\documentclass{article}

\PassOptionsToPackage{numbers}{natbib}
\usepackage[preprint]{neurips_2022}

\usepackage{times}
\usepackage[utf8]{inputenc}
\usepackage[T1]{fontenc}
\usepackage{url}
\usepackage{booktabs}
\usepackage{multirow}
\usepackage{amsfonts}
\usepackage{nicefrac}
\usepackage{microtype}
\usepackage{faktor}
\usepackage{hhline}

\usepackage{thm-restate}

\usepackage{tikz}
\usepackage{mathtools, stmaryrd}
\usepackage{tocloft}
\usepackage{enumitem}
\usepackage[hidelinks]{hyperref}
\usepackage{xcolor}
\hypersetup{
    colorlinks,
    linkcolor={red!50!black},
    citecolor={blue!50!black},
    urlcolor={blue!80!black}
}

\usepackage[algo2e, algoruled, boxed, vlined]{algorithm2e}
\usepackage{algorithm}
\usepackage{algpseudocode}

\SetKwInput{KwInput}{Input}
\SetKwInput{KwOutput}{Output}

\usepackage[noabbrev]{cleveref}

\title{Optimal first-order methods for convex functions \\ with a quadratic upper bound}

\author{
    Baptiste Goujaud\\
    CMAP, École Polytechnique, \\ Institut Polytechnique de Paris \\
    \texttt{baptiste.goujaud@polytechnique.edu}
    \And
    Adrien Taylor \\
    INRIA, École Normale Supérieure, \\ CNRS, PSL Research University, Paris \\
    \texttt{adrien.taylor@inria.fr}
    \And
    Aymeric Dieuleveut \\
    CMAP, École Polytechnique, \\ Institut Polytechnique de Paris \\ \texttt{aymeric.dieuleveut@polytechnique.edu}
}

\newtheorem{Th}{Theorem}[section]

\newtheorem{Def}[Th]{Definition}
\newtheorem{Assump}[Th]{Assumption}

\newtheorem{Rem}[Th]{Remark}
\newtheorem{?}[Th]{Problem}
\newtheorem{Ex}[Th]{Example}

\def\QG{\operatorname{QG}}
\def\RG{\operatorname{RG}}

\usepackage{amssymb}
\usepackage{pifont}
\newcommand{\cmark}{\textcolor{green!50!black}{\checkmark}}
\newcommand{\xmark}{\textcolor{red!50!black}{\times}}

\usepackage{float}
\usepackage{wrapfig}

\newcommand{\R}{\mathbb{R}}
\newcommand{\fs}{f_{\star}}

\begin{document}

\addtocontents{toc}{\protect\setcounter{tocdepth}{0}}

\maketitle

\begin{abstract}
    We analyze worst-case convergence guarantees of first-order optimization methods over a function class extending that of smooth and convex functions.
    This class contains convex functions that admit a simple quadratic upper bound.
    Its study is motivated by its stability under minor perturbations.
    We provide a thorough analysis of first-order methods, including worst-case convergence guarantees for several methods,
    and demonstrate that some of them achieve the optimal worst-case guarantee over the class.
    We support our analysis by numerical validation of worst-case guarantees using performance estimation problems.
    A few observations can be drawn from this analysis, particularly regarding the optimality (resp.~and adaptivity)
    of the heavy-ball method (resp.~heavy-ball with line-search).
    Finally, we show how our analysis can be leveraged to obtain convergence guarantees over more complex classes of functions.
    Overall, this study brings insights on the choice of function classes over which standard first-order methods have working worst-case guarantees.
\end{abstract}

\section{Introduction}\label{sec:introduction}
    In this paper, we consider the problem of minimizing a convex (closed proper) function
    \begin{equation}\label{eq:opt}
    \fs\triangleq\min_{x\in \R ^d} f(x),
    \end{equation}
    where $f: \R^d \to \R$ is assumed to have a non-empty set of global minimizers denoted by $\mathcal{X}_\star$ (which is necessarily convex).
    Convergence properties of first-order optimization are typically analyzed through worst-case analyses under the \emph{black-box} model~\citep{Book:NemirovskyYudin}.
    In this formalism, nontrivial guarantees are obtained by assuming the function to be minimized to satisfy certain regularity conditions.
    In particular, it is common to assume Lipschitz continuity of the gradients of $f$ (also often referred to as smoothness of $f$) as well as (strong) convexity of~$f$.
    First-order methods and their analyses for minimizing such functions in the black-box model occupied a great deal of attention, see, e.g.,~\citep{Book:NemirovskyYudin,Book:polyak1987,Nest03a}.

    Due to the practical success of first-order methods in large-scale applications,
    and particularly when aiming for only low to medium accuracy solutions (\citet{bottou2007tradeoffs} motivate this goal for machine learning),
    a few trends emerged in the first-order optimization literature.
    Among them, a particular focus concerned the question of understanding minimal working assumptions under which one could design efficient first-order methods.
    In other words, many authors looked for weaker/alternate versions to the standard smoothness and strong convexity-type assumptions,
    still allowing to obtain suitable working guarantees for standard first-order methods.

    \textbf{Relaxations of strong convexity-type assumptions.}
    Convexity alone is not sufficient to a priori guarantee ``fast'' convergence of usual first-order methods.
    On the other hand, strong convexity allows to obtain faster (geometric) rates but is a very strong condition.
    Therefore, many authors studied conditions in between convexity and strong convexity, aiming to obtain faster rates under relatively generic assumptions.
    In particular, different authors considered the restricted secant inequality~\citep{zhang2013gradient, guille2022gradient},
    the error bound~\citep{luo_error_1993}, {\L}ojasiewicz-type inequalities~\citep{polyak_gradient_1963}, and many more~\citep{NIPS2015_5718,kurdyka1998gradients,liu2014asynchronous,gong2014linear,necoara2016linear,hardt2016gradient,abbaszadehpeivasti2022conditions}.
    Relations between these assumptions were treated at length in~\citep{bolte2017error,zhang2017restricted}.
    Among those, one of the weakest relaxation is the so-called (lower) quadratic growth, see~\citep{bonnans1993second, ioffe1994sensitivity,anitescu2000degenerate}.
    Recently, those notions turned out to be useful, e.g.~for studying proximal gradient methods, see~\citep{cui2017quadratic, drusvyatskiy2016error, drusvyatskiy2015quadratic, zhang2020, zhang2017linear,chieu2021quadratic}.
    In the rest of the paper, we focus on (non strongly) convex functions.

    \textbf{Relaxations of smoothness-type assumptions.} Generalization of smoothness assumptions were less investigated in the literature.
    Still, a few such relaxations have emerged, including the {relative smoothness}, see~\citep{bauschke2017descent,lu2016relativelysmooth,dragomir2021optimal, hanzely2018accelerated}, restricted smoothness~\citep{agarwal2011fast} and restricted Lipschitz-continuous gradient~\citep{zhang2013gradient}.
    In this work, we consider instead the set of convex functions satisfying the (upper) \textit{quadratic growth} condition, as follows.

    \begin{Def}
        A function $f$ is $L$-quadratically upper bounded (denoted $L$-$\QG^+$) if for all $x \in \mathbb{R}^d$: \[f(x)-f_\star \leq \frac{L}{2} d(x, \mathcal{X}_\star)^2,\]
        where $d(x, \mathcal{X}_\star) = \min_{x_\star \in \mathcal{X}_\star} \|x-x_\star\|_2$.
        We denote the set of such functions by $\QG^+(L)$, and by $\QG^+$ when $L$ is left unspecified.
    \end{Def}

    This assumption is weaker than smoothness.
    First, any $L$-smooth function (i.e.\ with $L$-Lipschitz gradient) also belongs to $\QG^+(L)$.
    On the other hand,
    \begin{enumerate}[itemsep=1pt,topsep=0pt,leftmargin=*]
        \item some functions do belong to $\QG^+$ while not being smooth for any value of $L$.
        In particular $\QG^+$ contains all Lipschitz non-smooth convex functions which are twice differentiable at all their optimal points.
        Let us mention a few rules for obtaining (not necessarily smooth) $\QG^+$ functions:
        (i) any function that can be written as $x\mapsto h(x^\top M x)$, where $h$ is convex Lipschitz continuous and $M$ is a positive semidefinite matrix,
        or (ii)~$x\mapsto h(N(x))$, where $h$ is convex and smooth (or $\QG^+$) and $N$ is a norm (e.g. $N=\|\cdot\|_2$, $N=\|\cdot\|_1$ or $N=\|\cdot\|_\infty$).
        \item Some functions belong to $\QG^+(L_Q)$ while being smooth only for some $L_S \gg L_Q$
        (see e.g., Eq.(3) and Proposition 4.6 in~\citep{guille2021study}.)
        Consequently, even though the worst case convergence rate over the class of $\QG^+(L)$ functions cannot improve on the rate for $L$-smooth functions,
        it is possible for a given smooth function that the guarantee provided by the rate on the $\QG^+$ class is actually \textit{better} than the one resulting from the rate as a smooth function.
    \end{enumerate}

    More generally, the later example is related to \emph{condition continuity}, introduced by~\citet[Definition 4.9]{guille2021study}.
    \textit{Condition continuity} is a property of function classes, defined by the fact that a minor modification of the gradient of the function, at any point away from the optimum, cannot strongly affect the class parameter $L$.
    This is a desirable property for first-order methods for which the output typically continuously depends on the gradients of the functions minimized: if the function is slightly perturbed away from the optimum, the tuning and convergence guarantees of the algorithm should not be affected.
    The $\QG^+$ class satisfies condition continuity, while the class of smooth convex functions does not:
    a minor perturbation of an $L$-smooth function (thus $L$-$\QG^+$) can be $L_Q$-$\QG^+$ and $L_S$-smooth, with $L_S \gg L_Q$.
    This also motivates studying $\QG^+$.

    \textbf{Contributions and organization of the paper.}
    The rest of the paper is organized in two main sections.
    First, in Section~\ref{sec:algorithms_and_worst-case_analyses}, we analyze a few first-order methods, namely the subgradient method and the heavy-ball method with and without a line-search.
    We provide worst-case complexity bounds on the convergence rates  as well as corresponding lower complexity bounds.
    We also provide a lower complexity bound for minimizing convex functions in $\QG^+$ via first-order methods.
    Finally, we provide {interpolation/extensions} results for this class of problems, which allows exploring/deriving all previous results in a principled way (using performance estimation problems~\cite{drori2014performance,taylor2017smooth}).
    We summarize those results in \Cref{tab:summary_results}, together with precise references to the corresponding statements.
    Secondly, we review the main consequences of our analysis in \Cref{sec:discussion_take_away}.
    More specifically, we underline the facts that (a) the heavy-ball Algorithm~\ref{alg:ogm} and~\ref{alg:ogm_ls} are optimal on this class of functions, furthermore, (b) Algorithm~\ref{alg:ogm_ls} is  adaptive: it achieves the optimal convergence rate for both Lipschitz-continuous functions and $\QG^+$ functions without requiring knowledge of any class parameter.
    Then, we describe how our theory can be exploited for automatically obtaining convergence rates for different classes of functions.
    Lastly (in \Cref{apx:restart}), we discuss results that can be obtained when restricting the class to functions satisfying additional assumptions (a relaxation of strong-convexity).

    \begin{table*}
        {
        \caption{\label{tab:summary_results}
            Summary of the worst-case guarantees obtaining after $n$ iterations of a few different first-order methods on the class of convex $L$-$\QG^+$ functions, which are obtained in Section~\ref{sec:algorithms_and_worst-case_analyses} and \Cref{apx:subgrad}.
            Two main methods are studied, namely the (sub)gradient and the heavy-ball methods.
            Some guarantees concern the worst-case function value accuracy at the last iterate.
            The term ``average'' in the third column refers to the function value on the Polyak-Rupert averaged iterate.
            All the provided bounds are proportional to $R^2 = d(x_0, \mathcal{X}_\star)^2$, where $x_0$ denotes the starting point of the methods.
        }
        \begin{center}
            {\renewcommand{\arraystretch}{1.8}
            \resizebox{\linewidth}{!}{
             \begin{tabular}{@{}llrllrlr@{}}
                \specialrule{2pt}{1pt}{1pt}
                Method & \multicolumn{2}{l}{Step-sizes $(\gamma_t)_{0\le t\le n-1}$ } & Iterate & \multicolumn{2}{l}{Upper bound } & \multicolumn{2}{l}{Lower bound} \\
                \cmidrule{1-1}\cmidrule(l){2-3}\cmidrule(l){4-4}\cmidrule(l){5-6} \cmidrule(l){7-8}
                \multirow{3}{*}{Subgradient (Sec.~\ref{subsec:subgradient_method_on_qg+convex_functions})} & $\frac{1}{L}$ & (Alg.~\ref{alg:subgrad}) & Average & $\frac{L}{2}\frac{R^2}{n+1}$ & (Th.~\ref{thm:gd_average}) & $\frac{L}{2}\frac{R^2}{n+1}$ &(Rem.~\ref{rem:gd_average}) \\
                \cline{2-8}

                & $\gamma_t$ & (Alg.~\ref{alg:subgrad})& Last & $\times$ & & $\frac{LR^2}{2}L\gamma_{n-1}$ &(Th.~\ref{thm:non_convergence_gd}) \\
                \cline{2-8}

                & $\sim \frac{1}{2L\sqrt{t}}$& (Alg.~\ref{alg:gd_decreasing_step_sizes}) & Last & $\sim \frac{LR^2}{4\sqrt{n}}$ & (Conj.~\ref{conj:gd_sqrt}) & $ \sim \frac{LR^2}{4\sqrt{n}}$ & (Th.~\ref{thm:non_convergence_gd}) \\ \hline

                \multirow{2}{*}{Heavy-ball (Sec.~\ref{subsec:proposed_methods})} & $\frac{1}{L}\frac{1}{t+2}$ & (Alg.~\ref{alg:ogm}) & Last & $\frac{L}{2}\frac{R^2}{n+1}$ & (Th.~\ref{thm:hb_general}) & $\frac{L}{2}\frac{R^2}{n+1}$ & \\ \cline{2-8}

                & line-search &(Alg.~\ref{alg:ogm_ls}) & Last & $\frac{L}{2}\frac{R^2}{n+1}$ & (Th.~\ref{thm:hb_general}) & $\frac{L}{2}\frac{R^2}{n+1}$ & \\ \hline

                First-order (Sec.~\ref{subsec:first_order_lower_bound}) & Any& & Any & - &  & $\frac{L}{2}\frac{R^2}{n+1}$ &(Th.~\ref{thm:general_lower_bound}) \\

                \specialrule{2pt}{1pt}{1pt}\vspace{0em}
            \end{tabular}}
            \vspace{-1cm}
            }
        \end{center}}
    \end{table*}

    \textbf{Notation and background results.} For problem~\eqref{eq:opt}, the set $\mathcal{X}_\star$ of minimizers of $f$ is closed and convex ($f$ is proper closed and convex by assumption).
    Therefore, there exists a unique projection $\pi_{\mathcal{X_\star}}$ onto $\mathcal{X}_\star$, verifying:
    $\|x-\pi_{\mathcal{X_\star}}(x)\|_2 = d(x, \mathcal{X}_\star).$
    We use the classical $\partial f$ for denoting the subdifferential of the function $f$.
    Namely, the subdifferential of $f$ at $x\in\mathbb{R}^d$ is the set of all subgradients of $f$ at $x$:
    $\partial f(x) = \lbrace g | \forall y \in \mathbb{R}^d, f(y) \geq f(x) + \left< g | y - x \right> \rbrace.$
    Note that $0 \in \partial f(x) \Leftrightarrow x \in \mathcal{X}_\star$; moreover, if $f \in \QG^+$, then for all $x \in \mathcal{X}_\star$, $\partial f(x) = \lbrace 0 \rbrace$.

\section{\texorpdfstring{A few worst-case guarantees for minimizing $\QG^+$ convex functions}{A few worst-case guarantees for minimizing QG+ convex functions}} \label{sec:algorithms_and_worst-case_analyses}

    In this section, we provide the main technical results of this paper, summarized in \Cref{tab:summary_results}.
    In \Cref{subsec:subgradient_method_on_qg+convex_functions}, we study the behavior of a \emph{(sub)gradient} method on convex $\QG^+$ functions.
    A lower complexity bound on the convergence of any first-order method
    is provided in~\Cref{subsec:first_order_lower_bound}.
    In~\Cref{subsec:proposed_methods}, we introduce the heavy-ball method under consideration and prove its worst-case optimality.
    Finally, we discuss how interpolation conditions were used for obtaining these results in~\Cref{subsec:Interpol}.

    \subsection{\texorpdfstring{(Sub)gradient method on $\QG^+$ convex functions}{(Sub)gradient method on QG+ convex functions}}\label{subsec:subgradient_method_on_qg+convex_functions}
        
        \begin{wrapfigure}[5]{R}{0.4\textwidth}
            \vspace{-1.8cm}
            \begin{minipage}{0.4\textwidth}
                \begin{algorithm}[H]
                    \caption{Subgradient method}\label{alg:subgrad}
                    \KwInput{$x_0$, $(\gamma_t)_{0 \leq t \leq n}$}
                    \For{$k=1 \ldots n$}{
                        Query $g_{k-1} \in \partial f(x_{k-1})$;

                        $x_k \gets x_{k-1} - \gamma_{k-1} g_{k-1}$
                    }
                    \KwOutput{$(x_t)_{0 \leq t \leq n}$}
                \end{algorithm}
            \end{minipage}
        \end{wrapfigure}
        In this subsection, we consider \Cref{alg:subgrad}: the subgradient method, for $n$ iterations and a sequence of step-sizes $(\gamma_{t})_{0\le t\le n-1 }$.
        The following result provides a convergence guarantee for the averaged function value accuracy throughout the iterative procedure.

        \begin{restatable}{Th}{convergenceofsubgradientmethodinaverage}\textbf{\emph{(Convergence of \Cref{alg:subgrad} in average)}}
            \label{thm:gd_average}
            Let $f$ be an $L$-$\QG^+$ convex function.
            Applying (sub)gradient method on $f$ with step-size $\gamma \triangleq \frac{1}{L}$ leads to the following guarantee:
            \begin{equation}
                \frac{1}{n+1}\sum_{k=0}^n (f(x_k) - f_\star) \leq \frac{L}{2}\frac{1}{n+1} d(x_0, \mathcal{X}_\star)^2.
            \end{equation}
        \end{restatable}

        \noindent \textit{Sketch of proof.}
            The proof consists in proving that at each step
            $f(x_k) - f_\star \leq \frac{L}{2}d(x_k, \mathcal{X}_\star)^2 - \frac{L}{2}d(x_{k+1}, \mathcal{X}_\star)^2$
            and recognizing a telescopic sum on the right hand side.
            See Appendix~\ref{apx:subgrad}.
        $\hfill\blacksquare$

        By convexity of $f$, this result automatically implies a convergence guarantee for the Polyak-Ruppert~(PR) averaging~\citep{polyak1992acceleration, ruppert1988efficient}, $\bar x_n = \frac{1}{n+1}\sum_{k=0}^n x_k$ with the same convergence rate.
        For $L$-smooth convex functions, the same worst-case convergence rate is achieved by both the PR averaging and the last iterate.
        It is therefore natural to wonder if the subgradient method verifies the same convergence guarantee for the \emph{last iterate}, on $\QG^+$ convex functions.
        For $L$-smooth convex functions, \citet[Theorem 3.2]{drori2014performance} provide the following lower bound on the convergence of Algorithm~\ref{alg:subgrad}:
        for any $\gamma_k=\gamma \in \left[0, 2/L \right]$, there exists a $L$-smooth convex function and a starting point~$x_0$ s.t.
        \begin{equation}
            f(x_n) - f_\star \geq \max\left(\frac{L}{2}\frac{1}{1 + 2nL\gamma}, \frac{L}{2} \left( 1 - L\gamma \right)^{2n} \right)d(x_0, \mathcal{X}_\star)^2. \label{eq:gd_smooth}
        \end{equation}
        \citet[Theorem 3.1]{drori2014performance} also provide a corresponding worst-case guarantee of the form $f(x_n) - f_\star \leq \frac{L}{2}\frac{1}{1 + 2nL\gamma}d(x_0, \mathcal{X}_\star)^2$ for when $\gamma\in (0,\tfrac1L)$, which ensures convergence in function value accuracy with a constant step-size rule $\gamma_k=\gamma\in (0,\tfrac1L)$ at a rate $O(1/n)$.

        Here we provide a stricter lower bound for the convergence of the function value for the last iterate: contrary to what happens for smooth convex functions, subgradient methods with constant step-sizes cannot be guaranteed to converge on $\QG^+$ convex functions.

        \begin{restatable}{Th}{nonconvergenceofgdwithlowerboundedstepsizes}\textbf{\emph{(Lower bound for~\Cref{alg:subgrad} - final iterate).}}\label{thm:non_convergence_gd}
            \noindent For any sequence
            $(\gamma_i)_{0 \leq i \leq n-1}>0$ and any $\varepsilon>0$, there exists an $L$-$\QG^+$ convex function $f$ that verifies, after $n$ iterations of \Cref{alg:subgrad} with step-sizes $(\gamma_i)_{0 \leq i \leq n-1}$,
            \begin{equation} \label{eq:LB_QG}
                f(x_n) - f_\star \geq \frac{L}{2}L\gamma_{n-1} d(x_0, \mathcal{X}_\star)^2 - \varepsilon.
            \end{equation}
        \end{restatable}

        \noindent \textit{Sketch of proof.}
        The proof consists in finding a function defined on $\mathbb{R}^3$ such that all the iterates except the last one are very close to each other.
        As $\QG^+$ convex functions might not be differentiable, subgradients might vary very quickly.
        Therefore, the last iterate might be far away from the others, although all the previous iterates are clustered.
        A complete proof is provided in Appendix~\ref{apx:subgrad}.
        $\hfill\blacksquare$

        As a result, it is necessary to enforce $\gamma_n \rightarrow 0$ for ensuring convergence of Algorithm~\ref{alg:subgrad} on all problem instances.
        On the other hand, a similar lower bound to~\eqref{eq:gd_smooth} (using the same Huber function as that used for the class of smooth convex functions, see~\citep{drori2014performance}, or~\citep[Section 4]{taylor2017smooth}) and modifying $x_0$ to account for varying step-sizes (we use $x_0 = 1 + 2 \sum_{k=0}^{n-1} L\gamma_k$), one can obtain:
        \begin{equation}
            f(x_n) - f_\star \geq \frac{L}{2}\frac{1}{1 + 2 \sum_{k=0}^{n-1} L\gamma_k} d(x_0, \mathcal{X}_\star)^2.
            \label{eq:LB_smooth}
        \end{equation}
        Consequently the worst-case convergence is slower than $O(1 / n)$ as soon as $\gamma_k\rightarrow 0$.
        Overall, the convergence is provably worse for the last iterate over the $\QG^+(L)$-class than over the class of $L$-smooth convex functions, even though guarantees match for the PR-averaged iterate.
        Actually, the lower bound is at least the maximum of the RHSs of \Cref{eq:LB_QG,eq:LB_smooth}.
        In \Cref{apx:subgrad}, we introduce and analyze \Cref{alg:gd_decreasing_step_sizes}, that corresponds to \Cref{alg:subgrad} with a specific sequence of step-sizes such that \Cref{eq:LB_QG,eq:LB_smooth} are equal.
        This results in a decaying sequence of step-sizes scaling as $O(1/\sqrt{n})$.
        The next section is devoted to a lower complexity bound on the convergence in function accuracy for any black-box first-order method.

    \subsection{First-order lower bound}\label{subsec:first_order_lower_bound}

        The next theorem guarantees that no black-box first-order method can beat a $O(1 / n)$ worst-case guarantee in function values uniformly on the set of $\QG^+$ convex functions.

        \begin{center}
        \begin{restatable}{Th}{lowerboundoffirstorderalgorithm}\textbf{\emph{(Lower complexity bound)}}
            \label{thm:general_lower_bound}
            Let $n\in\mathbb{N}$.
            There exists some $d\in\mathbb{N}$ and some convex $L$-$\QG^+$ function $f$ of input space $\mathbb{R}^d$ such that: for any sequence $(x_k)_{0\le k\le n}$ satisfying $x_k - x_0 \in \text{span}\left\{g_0, g_1, g_2, \dots, g_{k-1}\right\}$ for all $k\leq n$ with $g_i \in \partial f(x_i)$ ($i=0,\ldots,k-1$), we have:
            \[f(x_n) - f_\star \geq \frac{L}{2}\frac{1}{n+1} d(x_0, \mathcal{X}_\star)^2.\]
        \end{restatable}
        \end{center}

        \noindent \textit{Sketch of proof.}
            The proof consists in noticing that the function $x \mapsto \frac{L}{2}\|x\|_\infty^2$ allows to explore one new dimension per step and that this new dimension is independent of the way the next point in the sequence is chosen.
            Therefore, choosing $d=n+1$ allows to ensure that there exists one \textit{unseen} dimension after $n$ iterations.
            This methodology is common to prove lower bounds in first-order optimization, see, e.g.,~\citep{nemirovskinotes1995,Nest03a,bubeck2014convex}.
            We refer to Appendix~\ref{apx:lower_bound_1} for the complete proof.
            The above result is also generalized in Appendix~\ref{apx:lower_bound_2} to account for \textit{any} sequence generated by a black-box first-order (possibly without the span assumption as in, e.g.,~\citep[Chapter 12]{nemirovskinotes1995} for quadratic minimization).
        $\hfill\blacksquare$

        One can conclude from \Cref{thm:general_lower_bound} that no black-box first-order method can enjoy a worst-case guarantee better than $f(x_n) - f_\star \leq \frac{L}{2}\frac{1}{n+1} d(x_0, \mathcal{X}_\star)^2$ uniformly on all $d\in\mathbb{N}$, all $f:\mathbb{R}^d\rightarrow\mathbb{R}$ that is convex and $L$-$\QG^+$ and all $x_0\in\mathbb{R}^d$.
        This entails that \Cref{alg:subgrad}, with constant step-size $1/L$ and PR averaging, is worst-case optimal for decreasing function values on the class of $\QG^+$ and convex functions.
        In the next section, we introduce two alternate methods that also achieve this optimal bound, this time without PR averaging.
        As we see in the sequel, those further developments allow achieving this optimal bound without explicitly using the knowledge of the constant $L$.
        
    \subsection{Two methods with optimal last iterate guarantee}\label{subsec:proposed_methods}

        \begin{wrapfigure}[8]{R}{0.5\textwidth}
            \begin{minipage}{0.5\textwidth}
                \vspace{-0.9cm}
                \SetInd{0.5em}{0.3em}
                \begin{algorithm}[H]
                    \caption{Heavy-ball method for $\QG^+$ convex}
                    \label{alg:ogm}
                    \KwInput{$x_0$, $L$}
                    \For{$k=1 \ldots n$}{
                    \text{Pick } $g_{k-1}$ from $\partial f(x_{k-1})$

                    $x_k \gets \frac{k}{k+1}x_{k-1} + \frac{1}{k+1}x_0 - \frac{1}{k+1}\sum_{i=0}^{k-1}\frac{1}{L} g_i$}
                    \KwOutput{$x_n$}
                \end{algorithm}
            \end{minipage}
        \end{wrapfigure}

        In this section, we introduce \Cref{alg:ogm} and \Cref{alg:ogm_ls} which both achieve the optimal convergence guarantee for the last iterate (see \Cref{thm:general_lower_bound}).
        The first of those two methods explicitly relies on the knowledge of the class parameter $L$ for performing its updates, whereas the second variant allows avoiding using any knowledge on $L$.
        Note that the update rule from \Cref{alg:ogm} can equivalently be expressed as:
        \begin{equation}
            x_k \gets x_{k-1} - \frac{1}{L}\frac{1}{k+1}g_{k-1} + \frac{k-1}{k+1}\left( x_{k-1} - x_{k-2} \right) \label{alg:hb_ogm}
        \end{equation}
        \begin{wrapfigure}[12]{R}{0.5\textwidth}
            \begin{minipage}{0.5\textwidth}
                \SetInd{0.5em}{0.3em}
                \vspace{-0.8cm}
                \begin{algorithm}[H]
                    \caption{Heavy-ball method with line-search for $\QG^+$ convex}
                    \label{alg:ogm_ls}
                    \KwInput{$x_0$, $v_0 \gets 0$}
                    \For{$k=1 \ldots n$}{
                        $y_k \gets \frac{k}{k+1}x_{k-1} + \frac{1}{k+1}x_0$

                        \text{Pick } $g_{k-1} \in \partial f(x_{k-1})$ such that $\left< g_{k-1} , v_{k-1} \right> = 0$.

                        $v_k \gets v_{k-1} + g_{k-1}$

                        $\alpha_k \gets \arg\min_{\alpha} f\left( y_k + \alpha v_k \right)$

                        $x_k \gets y_k + \alpha_k v_k$
                    }
                    \KwOutput{$x_n$}
                \end{algorithm}
            \end{minipage}
        \end{wrapfigure}
        where $g_{k-1} \in \partial f(x_{k-1})$.
        This formulation corresponds to the heavy-ball method, as defined in~\citep[Theorem 2]{ghadimi2014global} and for which  authors provided a $O(1/n)$ guarantee for $L$-smooth convex functions.
    
        \Cref{alg:ogm_ls} takes a similar form, but relies on an exact line-search procedure, avoiding to use any knowledge on $L$.
        Both methods share the same worst-case guarantee, matching the lower bound result from \Cref{thm:general_lower_bound}.
        The following theorem provides a necessary condition for an algorithm to share this same worst-case guarantee.

        \begin{center}
            \begin{restatable}{Th}{mainresult}\textbf{\emph{(Main result: sufficient condition for being worst-case optimal).}}
                \label{thm:main}
                Let $\mathcal{A}$ be an iterative first-order method that verifies, for all convex $\QG^+(L)$ function $f$, and starting points $x_0$,
                \begin{equation}
                    \Big\langle g_k , x_k - \Big[ \frac{k}{k+1}x_{k-1} + \frac{
                    1}{k+1}x_0 - \frac{1}{k+1}\sum_{i=0}^{k-1}\frac{1}{L} g_i \Big] \Big\rangle \leq 0. \label{eq:optimal_assumption}
                \end{equation}
                for some sequence $(g_i)_{i\in\mathbb{N}}$ of subgradients $g_i\in\partial f(x_i)$, and where $(x_i)_{i\in\mathbb{N}}$ are the iterates of $\mathcal{A}$.
                Then, the output $x_n$ of $\mathcal{A}$ achieves the worst-case guarantee:
                \begin{equation*}
                    f(x_n) - f_\star \leq \frac{L}{2}\frac{1}{n+1}d(x_0, \mathcal{X}_\star)^2.
                \end{equation*}
            \end{restatable}
        \end{center}

        \noindent \textit{Sketch of proof.}
            The proof is based on a Lyapunov analysis; for $k\geq 0$, we define the sequence
              $k(f(x_{k-1}) - f_\star) + \frac{L}{2} \| x_0 - \pi_{\mathcal{X_\star}}(x_0) - \sum_{i=0}^{k-1}\frac{1}{L} g_i \|^2$
            and show it is a decreasing.
            See Appendix~\ref{apx:main_result} for a complete and detailed proof.
        $\hfill\blacksquare$

        Inequality~\eqref{eq:optimal_assumption} is clearly satisfied for \Cref{alg:ogm} by ensuring the right hand side of the inner product being identically $0$.
        For \Cref{alg:ogm_ls}, the right hand side of the inner product in~\eqref{eq:optimal_assumption} is colinear to the search direction $\sum_{i=0}^{k-1} g_i$.
        First-order optimality conditions of the exact line-search procedure enforces the inner product in~\eqref{eq:optimal_assumption} to be identically $0$.
        The next corollary follows.

        \begin{center}
            \begin{restatable}{Cor}{optimalityofthetwoproposedalg} \label{cor:optimal}
            Let $n\in\mathbb{N}$, $d\in\mathbb{N}$, $f$ be a convex $\QG^+$ function, and $x_0\in\mathbb{R}^d$.
            Also let $x_n$ be the output of either Algorithm~\ref{alg:ogm} or Algorithm~\ref{alg:ogm_ls}, we have:
                $f(x_n) - f_\star \leq \frac{L}{2}\frac{1}{n+1}d(x_0, \mathcal{X}_\star)^2.$
            \end{restatable}
        \end{center}

        In the next section, we discuss how such worst-case analyses were obtained in a principled way, through so-called performance estimation problems (PEPs).
        An important ingredient to use this methodology is to develop \emph{interpolation} (a.k.a.~\emph{extension}) results for the  convex $\QG^+$ class.

    \subsection{\texorpdfstring{Extension/interpolation results for $\QG^+$ convex functions}{Extension-interpolation results for QG+ convex functions}}
    \label{subsec:Interpol}

        The problem of interpolating/extending within a class of functions can be stated as follows.
        Given a set of triplet $(x_i, g_i, f_i)_{i \in I}\subset \mathbb{R}^d\times\mathbb{R}^d\times \mathbb{R}$ for some $d\in\mathbb{N}$ and some index set $I$, the question of interest is that of recovering a function $f$ in a prescribed class of convex functions $\mathcal{F}$ satisfying
        \[ f_i=f(x_i) \text{ and } g_i\in\partial f(x_i) \text{ for all } i\in I.\]
        A similar problem, often referred to as convex integration, consists in finding such functions by only specifying some subgradients but no function values; see~\citep{Rock70}; for the case where $\mathcal{F}$ is the class of (closed and proper) convex functions, this problem was also treated at length in~\citep{lambert2004finite}.
        Motivated by applications to performance estimation problems (see below), this problem was studied in~\citep{taylor2017smooth} for the cases where $\mathcal{F}$ is the class of closed proper (possibly strongly) convex (possibly smooth) functions.
        In this case, it is possible to obtain simple necessary and sufficient conditions for the set $(x_i, g_i, f_i)_{i \in I}$ to be interpolable; we refer to those conditions as \emph{interpolation conditions}.
        Such conditions take the form of a set of inequalities on $(x_i, g_i, f_i)_{i \in I}$, and sometimes allow to conveniently deal with discrete versions of functions within a certain class $\mathcal{F}$ (for which we have interpolation conditions at our disposal).
        There exists a few classes of functions, typical for the analysis of first-order methods, for which such conditions exist, see, e.g.,~\citep[Theorem 3.3--3.6, Theorem 3.10]{taylor2017exact}.
        The next theorem provides interpolation conditions for the class of convex $\QG^+$ functions.

        \begin{center}
            \begin{restatable}{Th}{interpolationconditions}\textbf{\emph{(Interpolation conditions)}}
                \label{thm:interp}
                Let $(x_i, g_i, f_i)_{i \in I}$ a family of elements in $\mathbb{R}^d \times \mathbb{R}^d \times \mathbb{R}$.
                Set $I_\star$ the (assumed) non-empty subset of $I$ of the indices of elements $(x_i, g_i, f_i)$ verifying $g_i=0$.

                Then, there exists a $\QG^+(L)$ and convex function $f$ interpolating those points (i.e. such that $\forall i \in I, f(x_i) = f_i \text{ and } g_i \in \partial f(x_i)$) if and only if
                \begin{eqnarray}
                    \forall i \in I, \forall j \in I, && f_i \geq f_j + \left< g_j, x_i - x_j \right> \label{eq:interp_convexity} \\
                    \forall i \in I_\star, \forall j \in I, && f_i \geq f_j + \left< g_j, x_i - x_j \right> + \frac{1}{2L} \|g_j\|^2 .\label{eq:interp_convexity_qg}
                \end{eqnarray}
            \end{restatable}
        \end{center}

        \noindent \textit{Sketch of proof.}
            The proof is derived in two steps.
            First we notice that~\eqref{eq:interp_convexity} corresponds to the convexity of the function, and we prove~\eqref{eq:interp_convexity_qg} combining the $2$ inequalities respectively corresponding to convexity and $L$-$\QG^+$ assumptions.
            Reciprocally, we explicitly build a $L$-$\QG^+$ convex function from~\eqref{eq:interp_convexity} and~\eqref{eq:interp_convexity_qg}.
            See Appendix~\ref{apx:interpolation_conditions} for a detailed proof.
        $\hfill\blacksquare$

        \textbf{Application to Performance estimation problems (PEPs).} PEPs were introduced by~\citet{drori2014performance} for developing new analyses of first-order methods; see also~\citep{ drori2014contributions, kim2016optimized} for the first works on this topic.
        PEPs were later formalized using the concept of \emph{convex interpolation} by~\citep{taylor2017smooth, taylor2017exact}.
        PEPs formulate the search for worst-case guarantees as infinite dimensional optimization problems over the considered class of functions, e.g.,
        \[
        \text{Worst case}\left(\text{Alg.}~\ref{alg:subgrad}, n=1, \gamma=\tfrac{1}{L}, \text{ convex } \QG^+(L)  \right) \triangleq \max_{\footnotesize\shortstack{$f\in \QG^+(L)$ \text{ convex } \\ $x_1 \in x_0 -\frac{1}{L} \partial f(x_0)$} } \frac{f(x_1) - f_\star}{ d(x_0,\mathcal X_\star)^2},
        \]
        for the case of 1-step subgradient descent.
        In order to numerically solve those problems, it is  needed to transform them into a finite dimensional problem.
        To that end, interpolation conditions play a crucial role, by allowing to reduce the optimization over the (infinite dimensional) class of functions to an optimization over a constrained set of vectors, thus finite dimensional problem.

        Consequently, \Cref{thm:interp} allows us to use the PEP framework to study the class of $\QG^+(L)$ convex functions.
        Therefore, we obtained and verified all the results of this work thanks to the different programming tools that have recently been developed for easing the access to the PEP framework (see the packages~\citep{taylor2017performance, goujaud2022pepit}).
        For each algorithm studied in this paper, PEPs also provided the associated tight bound.
        Finally, PEPs also guided us to prove the lower bound of Theorem~\ref{thm:general_lower_bound} thanks to the study of the greedy first-order method (\emph{GFOM}) from~\citep{drori2020efficient}.

    \vspace{-0.5em}

\section{Discussion and concluding remarks}
\label{sec:discussion_take_away}

    \vspace{-0.5em}
    In this section, we discuss a few takeaways of the results.
    The messages of this section include optimality and adaptivity results for heavy-ball with a line-search, a discussion on the applicability of this functional class beyond its simple use, as well as a few words on the limitations of considering this simple class of convex functions.

    \subsection{Optimality of HB algorithm}
    \label{subsec:optimality-of-hb-algorithm}

        An ever-recurring question in the field of optimization is the convergence of HB methods on smooth and strongly convex functions.
        On the one hand, HB is optimal (in the sense that it achieves the optimal worst-case guarantee) for convex quadratic objectives and can then be seen as a variant of Chebyshev iterative method~\citep{flanders1950numerical,lanczos1952solution, young1953richardson}.
        Even with a simple (constant) tuning of the step-size and momentum parameters, it is optimal (i.e., achieving rates $O(1/n^2)$ for $L$-smooth convex quadratics, and $(1-O(\sqrt{\mu/L}))^n$ if the problem is also $\mu$-strongly convex.).

        On the other hand, the  method does not generalize well to (non quadratic) smooth strongly convex functions: \citet[Figure 7]{lessard2016analysis} built a function for which the heavy-ball method, tuned with the same dependence to $L$ and $\mu$ than for the quadratic case, fails to converge.
        More generally, while other sets of hyper-parameters allow to obtain convergence for the class of smooth and (strongly) convex~\citep{ghadimi2014global}, heavy-ball was never showed to accelerate w.r.t.~the gradient descent method, and the possibility to obtain such an acceleration remains an open question to the best of our knowledge.
        Simultaneously, Nesterov's accelerated gradient method~\citep{nesterov1983method} does achieve such an acceleration.

        In summary, while for quadratic convex problems it has the optimal worst-case guarantee, heavy-ball is believed not to satisfy this property for smooth and (strongly) convex functions.
        Interestingly, \Cref{cor:optimal} shows that heavy-ball is optimal (with rate $O(1/n)$) over the (larger) class $\QG^+$ convex.

        This observation questions the existence of intermediary classes (smaller than $\QG^+$ convex and containing quadratic functions), over which heavy-ball would achieve a $O(1/n^2)$ convergence guarantee.
        In the next section, we discuss the adaptivity of the method.

    \subsection{\texorpdfstring{Adaptivity of HB line-search algorithm 3}{Adaptivity of HB line-search algorithm~\ref{alg:ogm_ls}}}
    \label{subsec:adaptivity_hb_line-search}

        The search for adaptive and parameter free methods is a major challenge in optimization as the regularity of the function (both in terms of class and class-parameter $L$) is often unknown.
        In the rest of the section, we discuss the fact that~\Cref{alg:ogm_ls} provides a parameter-free and adaptive to $\QG^+(L)$ and $M$-Lipshitz functions method.
        Those results are summarized in \Cref{tab:optimality_summary}.

        \begin{table*}
            {
            \caption{\label{tab:optimality_summary}
                Optimality of the proposed methods over the set of $\QG^+$ convex functions and $M$-Lipschitz convex functions.
                ELS: Exact Line-Search. $\cmark$ indicates optimality among the class and $\xmark$ the contrary.
                All counter examples are given in App.~\ref{apx:tab}. $\ ^{\dagger}$:~constants resulting in optimal convergence rates depend on the class, thus for example heavy-ball with step-size $\frac{\text{constant}}{(t+2)}$ is not adaptive as it does not achieve the optimal rate for both classes with the same constant. $\ ^\ddagger$: up to a $\log$ factor.
                } \vspace{-0.3cm}
            \begin{center}
                {\renewcommand{\arraystretch}{1.8}
                \resizebox{\linewidth}{!}{
                 \begin{tabular}{@{}lllcrcrc@{}}
                    \specialrule{2pt}{1pt}{1pt}
                    \multicolumn{3}{c}{Method} & \multicolumn{4}{c}{Function class} & Parameter free \\
                     \cmidrule{1-3}\cmidrule(l){4-7}\cmidrule(l){8-8}
                    Algorithm & Step-sizes $(\gamma_t)_{0\le t\le n-1}$  & Iterate & \multicolumn{2}{c}{$\QG^+(L)$ convex} & \multicolumn{2}{c}{$M$-Lipschitz convex}   \\
                    \cmidrule(l){4-5}\cmidrule(l){6-7}
                    Subgradient (Alg.~\ref{alg:subgrad}) & $\text{constant}^{\dagger}$ & Average &  $\cmark$ & (Thm.~\ref{thm:gd_average}) & $ \xmark$ & (Thm.~\ref{thm:subgrad_constant_lower_lip}) & $\xmark$ \\
                    Subgradient (Alg.~\ref{alg:gd_decreasing_step_sizes}) & $ \text{constant}^{\dagger}/\sqrt{t}$ & Average & $\xmark$ & \eqref{eq:LB_smooth} & $\cmark^\ddagger$ &~\citep[Sec. 3.2.3]{Nest03a} & $\xmark$ \\
                    Subgradient (Alg.~\ref{alg:subgrad_els}) & ELS & Average &  $\xmark$ & (Thm.~\ref{thm:gd_els_lower_bound}) & $\xmark$ & (Thm.~\ref{thm:gd_els_lower_bound}) & $\cmark$ \\
                    Subgradient (Alg.~\ref{alg:subgrad_els}) & ELS & Last &  $\xmark$ & (Thm.~\ref{thm:gd_els_lower_bound}) & $\xmark$ & (Thm.~\ref{thm:gd_els_lower_bound}) & $\cmark$ \\
                    \midrule
                    Heavy-ball (Alg.~\ref{alg:ogm}) & $\text{constant}^{\dagger}/(t+2)$ & Last & $\cmark$ & (Cor.~\ref{cor:optimal}) & $\cmark$ &~\citep[][, Cor. 3]{drori2020efficient} & $\xmark$ \\
                    Heavy-ball (Alg.~\ref{alg:ogm_ls}) & ELS & Last & $\cmark$ & (Cor.~\ref{cor:optimal}) & $\cmark$ &~\citep[][, Cor. 4]{drori2020efficient} & $\cmark$ \\
                    \specialrule{2pt}{1pt}{1pt}\vspace{0em}
                \end{tabular}}
                \vspace{-1cm}
                }
            \end{center}}
        \end{table*}

        \textbf{Non convergence of the line-search version of \Cref{alg:subgrad}.}
        While the line-search version of Algorithm~\ref{alg:ogm} (i.e. Algorithm~\ref{alg:ogm_ls}) allows to get rid of the knowledge of $L$ without degrading the convergence guarantee for the latest iterate, this is not the case for subgradient method.
        In fact subgradient with line-search \Cref{alg:subgrad_els} does not converge, even when looking at the PR averaged iterate (see \Cref{thm:gd_els_lower_bound}).
        It thus destroys the  guarantee given in Theorem~\ref{thm:gd_average} for Algorithm~\ref{alg:subgrad} with constant step-size $1/L$.

        \textbf{\Cref{alg:subgrad} is not adaptive to the class of Lipschitz functions.} While \Cref{alg:subgrad} with averaging and constant step-size $1/L$ is optimal for the class $\QG^+(L)$ convex, for any sequence of constant step-sizes, neither the average nor the last iterate of \Cref{alg:subgrad} converge over the class of $M$-Lipschitz functions (see \Cref{thm:subgrad_constant_lower_lip}).
        On the other hand, to obtain a nearly (up to a log factor) optimal convergence rate $O(\log(n) / \sqrt{n})$, one can use $\gamma_t = \text{constant}/\sqrt{t}$ for $t\in \{1, \dots, n\}$, but such a sequence degrades the converge for $\QG^+$ as proved by the lower bound~\Cref{thm:non_convergence_gd}.

        \textbf{Adaptivity of HB line-search (Algorithm~\ref{alg:ogm_ls}).}
        While Algorithm~\ref{alg:ogm_ls} requires to perform exact line-search steps, its first advantage over Algorithm~\ref{alg:ogm} is not requiring the knowledge of the class parameter~$L$. \Cref{alg:ogm_ls} is thus  adaptive to the class parameter $L$.
        This is analogous of \emph{SSEP-based subgradient method} presented in \citet[Corollary 3]{drori2020efficient} and its line-search version~\citep[See][Corollary 4]{drori2020efficient} that are both optimal on the class of Lipschitz continuous and convex functions.
        The first one requires the knowledge of the parameter class $M$, the distance of the starting point to optimum $d(x_0, \mathcal{X}_\star)$ as well as the number of performed steps, while the second one replaces this knowledge by exact line-search steps.
        Moreover, we note that \emph{SSEP-based subgradient method} and Algorithm~\ref{alg:ogm} are very similar to each other.
        Indeed, their respective update steps can be written as
        $x_k \gets \frac{k}{k+1}x_{k-1} + \frac{1}{k+1}x_0 - \frac{d(x_0, \mathcal{X}_\star)}{M\sqrt{N+1}} \frac{1}{k+1}\sum_{i=0}^{k-1} g_i$
        and $x_k \gets \frac{k}{k+1}x_{k-1} + \frac{1}{k+1}x_0 - \frac{1}{L}\frac{1}{k+1}\sum_{i=0}^{k-1} g_i.$

        Remarkably, only the two constants in front of $\frac{1}{k+1}\sum_{i=0}^{k-1} g_i$ differ between the  two methods.
        Thus, their two line-search versions are identical.
        We conclude from Corollary~\ref{cor:optimal} and \citet[Corollary 4]{drori2020efficient} that Algorithm~\ref{alg:ogm_ls} is optimal on the classes of Lipschitz convex functions \emph{and} on the classes of $\QG^+$ convex functions.
        One can run this algorithm without knowing the type of conditions that are verified by the objective function and still benefit from a guarantee of optimality.
        This is a significant argument in favor of \Cref{alg:ogm_ls}.

    \subsection{Leveraging our analysis to obtain convergence bounds on other classes}
    \label{subsec:leveraging-our-analysis-to-obtain-convergence-bounds-on-other-classes}
    
        One of the major limitations of the $\QG^+$ class is that all functions within the class must be twice differentiable at the set of optimal points.
        This typically excludes some Lipschitz functions, as for example $x \mapsto \|x\|_1$ or problems involving Lasso regularization.
        In this section, we show that our \Cref{sec:algorithms_and_worst-case_analyses} can be leveraged to automatically obtain rates on more complete classes of functions, that combine the limitations of the Lipschitz convex and $\QG^+$ convex classes.\\
        To introduce these classes, we denote $h$ a generic function defined on $\mathbb{R}^+$ and verifying: $h(0)=0$, $h$ is strictly increasing and $h$ is concave.
        We consider the class of functions with $h$ \textit{relative growth}.
        \begin{Def}
            A function $f$ is $h$-relatively upper bounded (denoted $h$-$\RG^+$) if for all $x \in \mathbb{R}^d$: \[f(x)-f_\star \leq h\left( d(x, \mathcal{X}_\star)^2 \right).\]
            We denote the set of such functions by $\RG^+(h)$, and by $\RG^+$ when $h$ is left unspecified.
            \label{def:h_rg}
        \end{Def}
        These classes enable to cover the $\QG^+$ classes, the Lipschitz classes, and more.
        \begin{Rem}(Examples)
            \label{rem:h_linear}
            These three examples of functions $h$ satisfy the required assumptions:
            \begin{enumerate}[noitemsep, topsep=1pt,leftmargin=*]
                \item When $h$ is the linear function $h: z \mapsto \frac{Lz}{2}$, then simply $\RG^+(h) = \QG^+(L)$.
                \item When $h$ is the function $h: z \mapsto M\sqrt{z}$, then simply $\RG^+(h) = \left\{ M\text{-Lipschitz continuous}\right\}$.
                \item When $h: z \mapsto M\sqrt{z} + \frac{Lz}{2}$ a broader class  containing the limitations of both previous ones.
            \end{enumerate}
        \end{Rem}
        In the following, we consider the class of $h$-$\RG^+$ and convex functions.
        We then propose Algorithm~\ref{alg:hb_general}, which consists in applying \Cref{alg:ogm} (with the update written as in \Cref{alg:hb_ogm}) to the $\QG^+$ convex function $h^{-1} \circ (f - f_\star)$.
        We obtain the following convergence rate.
    
        \begin{figure}[h!]
            \begin{minipage}{0.49\linewidth}
                \SetInd{0.5em}{0.3em}
                \begin{algorithm}[H]
                    \caption{Heavy-ball method for $h$-$\RG^+$ convex functions}
                    \label{alg:hb_general}
                    \KwInput{$x_0$, $h$, $f_\star$}
                    \For{$k=1 \ldots n$}{
                    Choose ~ $g_{k-1}$ from $\partial f(x_{k-1})$

                    $\gamma_{k-1}= \frac{1}{2(k+1)}\frac{1}{h' \circ h^{-1}\left(f(x_{k-1}) - f_\star\right)} $

                    $x_k\! \gets\! x_{k-1} - \gamma_{k-1} g_{k-1} + \frac{k-1}{k+1}\left( x_{k-1} - x_{k-2} \right)$
                    }
                    \KwOutput{$x_n$}
                \end{algorithm}
            \end{minipage}
            \hspace{0.02\linewidth}
            \begin{minipage}{0.49\linewidth}
                \SetInd{0.5em}{0.3em}
                \begin{algorithm}[H]
                    \caption{Heavy-ball method for Lipschitz continuous convex functions}
                    \label{alg:hb_lip_from_qg}
                    \KwInput{$x_0$, $M$, $f_\star$}
                    \For{$k=1 \ldots n$}{
                    Choose ~ $g_{k-1}$ from $\partial f(x_{k-1})$

                    $\gamma_{k-1}=\frac{1}{k+1}\frac{f(x_{k-1}) - f_\star}{M^2}$

                    $x_k\! \gets\! x_{k-1} - \gamma_{k-1} g_{k-1} + \frac{k-1}{k+1}\left( x_{k-1} - x_{k-2} \right) $
                    }
                    \KwOutput{$x_n$}
                \end{algorithm}
            \end{minipage}
        \end{figure}
    
        \begin{center}
            \begin{restatable}{Th}{hbgeneral}\textbf{\emph{\label{thm:hb_general}}}
                Let $f$ an $h-\RG^+$ convex function, and $x_0$ any starting point.
                Then Algorithm~\ref{alg:hb_general} verifies
                    $f(x_n) - f_\star \leq h\left(\frac{d(x_0, \mathcal{X}_\star)^2}{n+1} \right).$
                In the examples of~\Cref{rem:h_linear}, this gives:
                \begin{enumerate}[noitemsep, topsep=1pt,leftmargin=*]
                    \item When $h$ is the linear function $h: z \mapsto \frac{Lz}{2}$ ($f$ is $L$-$\QG^+$ convex), then $f(x_n) - f_\star \leq \frac{L}{2}\frac{d(x_0, \mathcal{X}_\star)^2}{n+1}$.
                    \item When $h$ is the function $h: z \mapsto M\sqrt{z}$ ($f$ is $M$-Lip. convex), then $f(x_n) - f_\star \leq M\frac{d(x_0, \mathcal{X}_\star)}{\sqrt{n+1}}$.
                    \item When $h$ is the function $h: z \mapsto M\sqrt{z} + \frac{Lz}{2}$, then $f(x_n) - f_\star \leq M\frac{d(x_0, \mathcal{X}_\star)}{\sqrt{n+1}} + \frac{L}{2}\frac{d(x_0, \mathcal{X}_\star)^2}{n+1}$.
                \end{enumerate}
            \end{restatable}
        \end{center}

        \noindent \textit{Sketch of proof.}
        The proof consists in inverting $h$ and applying one of the results on the $\QG^+$ convex functions $h^{-1}\left( f - f_\star \right)$.
        For a detailed proof, see Appendix~\ref{apx:upper_assumption}.
        $\hfill\blacksquare$

        We also observe that for $h: z \mapsto \frac{Lz}{2}$, Algorithm~$\ref{alg:hb_general}$ is exactly Algorithm~$\ref{alg:ogm}$ and we recover the worst-case guarantee provided in Theorem~\ref{thm:hb_general}.
        Similarly, when considering $h: z \mapsto M\sqrt{z}$, Algorithm~$\ref{alg:hb_general}$ is written as \Cref{alg:hb_lip_from_qg} and we recover the worst-case guarantee provided in~\citep[Cor. 3]{drori2020efficient}.
        However, in this case the algorithm is different: the latest requires the knowledge of $d(x_0, \mathcal{X}_\star)$, while \Cref{alg:hb_lip_from_qg} requires the knowledge of $f_\star$.
        In this sense, the two proposed methods are complementary.

        Finally, we emphasize the fact that in practice, a lot of machine learning models requires to minimize non-smooth functions that are neither Lipschitz continuous nor $\QG^+$ (e.g. least-square regressions with lasso penalization of $\text{TV-L}_2$ model widely used in computer vision).
        These functions can be tackled using the flexible function $h(z) = M\sqrt{z} + \frac{L}{2}z$
        Then, applying \Cref{thm:hb_general} leads to the guarantee given in the third point of \Cref{thm:hb_general} obtained with \Cref{alg:hb_general} where $\gamma_{k-1} \gets \frac{1}{k+1}\frac{1}{L} \Big[ 1 - 1/\sqrt{1 + \frac{2L}{M^2}(f(x_{k-1}) - f_\star)} \Big]$.
    
        \textbf{Extension with additional constraint.} In \Cref{apx:restart}, we provide geometric convergence guarantees when the functions are also assumed satisfy a relaxation of strong convexity, $\QG^-$.

        \textbf{List of potential applications.}
        The $\QG^+$ class, and its extension to $h-\RG^+$ offer a flexibility that allows to tackle several non-smooth machine learning problems, including for example RELU activation, $L_1$ or TV regularization.
        As an example, the classical TV-$L_2$ denoising problem, which combines a term with quadratic growth with a non-smooth (even at the optimum) term, is neither Lipschitz nor smooth, but belongs to our class $h-\RG^+$, with $h: z \mapsto M \sqrt{z}+\frac{L z}{2}$ (\Cref{rem:h_linear}-3).

        \textbf{Limitations.}
        As specific methods have often been designed for those applications, our approach does not bring a systematic improvement.
        Yet, we believe it paves the way for adaptive methods that could do so.
        Remark that our Algorithm~$\ref{alg:ogm_ls}$ is adaptive and optimal for both the class of $\QG^+(L)$ functions, and the class of $M$-Lipschitz continuous ones.
        Extending our analysis to obtain an adaptive algorithm for $h-\RG^+$, which is left as an open direction, would allow to efficiently tackle TV-$L^2$ type of problems in a parameter-free way.

        \paragraph{Conclusion.}
        In this paper, we thoroughly analyze the class of convex $\QG^+$ functions.
        This function class relaxes the smoothness assumption and is motivated by the fact that $\QG^+$ satisfies \textit{condition continuity}, a desirable property for analyzing first-order methods.
        We analyze several such methods, and provide tight worst-case guarantees for them.
        Three methods achieve the optimal convergence rate over the class.
        Our analysis is supported by the derivation of \textit{interpolation conditions} allowing to verify all the results numerically.
        In particular, we observe that a heavy-ball algorithm results in acceleration (w.r.t.~the subgradient method), attaining the lower complexity bound for this class, a surprising result with respect to the smooth case.
        Moreover, using line-search, we obtain a parameter-free algorithm which is adaptive to the class parameter in $\QG^+$, and to the function class, as it also achieves the optimal rate for Lipschitz functions, a strongly desirable property in practice.
        Finally, we leverage our results to obtain convergence bounds for more complex classes of functions, combining the difficulties of the $\QG^+$ and the Lipschitz classes.
        Overall, this work participates to the trend of questioning the relevance of the most classical assumptions used in the analysis of first-order optimization methods.
        While our results are not intended to provide a definitive answer to this question, it goes one step further by providing an in-depth analysis for a more stable class of functions.
        Providing a similar analysis while relaxing convexity is a major open challenge.

\begin{ack}
    The work of B. Goujaud and A. Dieuleveut is partially supported by ANR-19-CHIA-0002-01/chaire SCAI, and Hi!Paris.
    A. Taylor acknowledges support from the European Research Council (grant SEQUOIA 724063).
    This work was partly funded by the French government under management
    of Agence Nationale de la Recherche as part of the ``Investissements d’avenir'' program,
    reference ANR-19-P3IA-0001 (PRAIRIE 3IA Institute).
\end{ack}

\bibliographystyle{plainnat}
\bibliography{references}

\clearpage
\appendix

    \section*{Organisation of the appendix}
    
        This appendix contains the proofs of the theorems stated in the main core of the paper.
        We also state a conjecture and bring some evidence about its statement.
        This appendix also contains discussions and extended results.
        
        Appendix~\ref{apx:subgrad} details the results on the subgradient method.
        Appendix~\ref{apx:gd_proof_th1_upper_bound} contains the proof of \Cref{thm:gd_average},
        Appendix~\ref{apx:gd_lower_bound} contains the proof of \Cref{thm:non_convergence_gd}
        and Appendix~\ref{apx:gd_conjecture} contains a conjecture that does not appear in the main core of the paper.
        This appendix also contains some evidence supporting this conjecture.
        
        Appendix~\ref{apx:lower_bound} contains the proofs for lower bounds on the class $\QG^+$ convex.
        Appendix~\ref{apx:lower_bound_1} contains the proofs of \Cref{thm:general_lower_bound}
        stating the lower bound under the classical assumption
        that the difference between the iterates lies into the span of observed gradients.
        Appendix~\ref{apx:lower_bound_2} extends the latter results without the aforementioned assumption.
        
        Appendix~\ref{apx:main_result} contains the proof of Theorem~\ref{thm:main},
        the main result of the paper, stating that all first order algorithm verifying a given identity,
        also enjoys an upper bound guarantee.
        
        Appendix~\ref{apx:tab} contains the proofs of all the claims that figure in \Cref{tab:optimality_summary}
        that are not already made elsewhere in this work or in others.
        
        Appendix~\ref{apx:interpolation_conditions} contains the proof of \Cref{thm:interp},
        essential to use the PEP framework.
        
        Appendix~\ref{apx:upper_assumption} contains all the proofs
        and discussions related to the extended class of the $\RG^+$ convex functions
        
        Finally, Appendix~\ref{apx:restart} contains linear convergence result
        under an additional assumption similar to the classical quadratic growth assumption.
        This result is not presented in the main core in the paper,
        since it is a bit out of the scope of the main message.
        However, we thought it was worth mentioning it here.
    
    \hypersetup{linkcolor = black}
    \setlength\cftparskip{2pt}
    \setlength\cftbeforesecskip{2pt}
    \setlength\cftaftertoctitleskip{3pt}
    \addtocontents{toc}{\protect\setcounter{tocdepth}{2}}
    \setcounter{tocdepth}{1}
    \tableofcontents
    \hypersetup{linkcolor=blue}
    
    \section{\texorpdfstring{(Sub)gradient method on $\QG^+$-convex functions}{(Sub)-gradient method on QG+ convex functions}}
\label{apx:subgrad}

    In this appendix, we provide the proof of Theorems~\ref{thm:gd_average} and~\ref{thm:non_convergence_gd} stating respectively an upper bound result on the subgradient method with fixed step-size $1 / L$ on the Polyak-Rupert averaged iterate, and a lower bound result on the subgradient method on the last iterate.
    Finally, based on this lower bound, we suggest a specific tuning of the subgradient method for $\QG^+$ convex functions.
    A conjecture is formulated on the worst-case bound achieved by this method with the prescribed tuning, as well as evidence obtained through the PEP framework.
    
    \subsection{Convergence of subgradient method with fixed step-size at Polyak-Rupert averaged iterate}
    \label{apx:gd_proof_th1_upper_bound}

        In section~\ref{subsec:subgradient_method_on_qg+convex_functions}, we state the following theorem about a worst-case upper bound of \Cref{alg:subgrad} on the class of $\QG^+$-convex functions.
        In this section, we provide the proof of this theorem.
        
        \convergenceofsubgradientmethodinaverage*
    
        \noindent \textit{Proof.}
            Let $k \in [0, n]$.
            We have
            \begin{eqnarray*}
                d(x_{k+1}, \mathcal{X}_\star)^2 & = & \|x_{k+1} - \pi_{\mathcal{X}_\star}(x_{k+1})\|^2 \\
                & \leq & \|x_{k+1} - \pi_{\mathcal{X}_\star}(x_{k})\|^2 \\
                & = & \|x_k - \gamma g_k -  \pi_{\mathcal{X}_\star}(x_{k})\|^2, \qquad \text{ with } g_k \in \partial f(x_k) \\
                & = & \|x_k - \pi_{\mathcal{X}_\star}(x_{k})\|^2 - 2\gamma \left< x_k - \pi_{\mathcal{X}_\star}(x_{k}) , g_k \right> + \gamma^2 \|g_k\|^2 \\
                & \overset{Eq.~\eqref{eq:interp_convexity_qg}}{\leq} & \|x_k - \pi_{\mathcal{X}_\star}(x_{k})\|^2 - 2\gamma \left(f(x_k) - f_\star + \frac{1}{2L}\|g_k\|^2\right) + \gamma^2 \|g_k\|^2 \\
                & = & d(x_k, \mathcal{X}_\star)^2 - 2\gamma \left(f(x_k) - f_\star \right) - \gamma \left( \frac{1}{L} - \gamma \right) \|g_k\|^2 \\
                & \overset{\gamma=\frac{1}{L}}{=} & d(x_k, \mathcal{X}_\star)^2 - \frac{2}{L} \left(f(x_k) - f_\star \right)
            \end{eqnarray*}
            
            By reordering the terms and summing over $k$:
            \begin{equation}
                \sum_{k=0}^n (f(x_k) - f_\star) \leq \frac{L}{2}d(x_0, \mathcal{X}_\star)^2
            \end{equation}
            
            which leads to the desired results.
            
        $\hfill\blacksquare$
    
        \begin{Rem}
            \label{rem:gd_average}
            From Theorem~\ref{thm:gd_average}, we conclude
            \begin{eqnarray*}
                \min_{0 \leq k \leq n} f(x_k) - f_\star & \leq & \frac{L}{2}\frac{1}{n+1} d(x_0, \mathcal{X}_\star)^2 \\
                f\left( \frac{1}{n+1}\sum_{k=0}^n x_k \right) - f_\star & \overset{\text{(by convexity)}}{\leq} & \frac{L}{2}\frac{1}{n+1} d(x_0, \mathcal{X}_\star)^2.
            \end{eqnarray*}
        \end{Rem}
    
        \begin{Rem}
            Note that this bound is tight not only for $\QG^+$ convex functions, but also for smooth convex functions.
            
            Indeed, we consider the real Huber function defined as
            \begin{equation}
                \label{eq:huber}
                f(x) = \left\{
                \begin{array}{cc}
                    \frac{L}{2}x^2 & \text{if } |x| \leq 1 \\
                    L|x| - \frac{L}{2} & \text{if } |x| \geq 1
                \end{array}
                \right.
            \end{equation}
            
            This function is $L$-smooth convex and often used to find lower bounds~\citep[See e.g.][]{drori2014performance,taylor2017smooth,kim2016optimized}.
            Moreover, starting from $x_0 = n + 1$, GD with $\gamma = \frac{1}{L}$ leads exactly to $x_k = n+1-k$ for all $k\leq n$, hence $f(x_k) - f_\star = L \left(n+\frac{1}{2}-k\right)$ and $\sum_{k=0}^n (f(x_k) - f_\star) = L \sum_{k=0}^n n+\frac{1}{2}-k = \frac{L}{2}(n+1)^2 = \frac{L}{2}d(x_0, \mathcal{X}_\star)^2$.
        \end{Rem}
    
    \subsection{Convergence limitation of the subgradient method in last iterate}
    \label{apx:gd_lower_bound}

        In this section, we prove Theorem~\ref{thm:non_convergence_gd} stating a lower bound guarantee on the convergence of the subgradient method.
        This proof is by far the most technical of this paper due to the amount of newly introduced notations.
        
        \nonconvergenceofgdwithlowerboundedstepsizes*
        
        \noindent \textit{Proof.}
            Let $\eta > 0$.
            
            We introduce the following notations:
            \begin{itemize}
                \item $\delta \triangleq \left(\frac{\eta\sqrt{3}}{1+L\gamma_{n-2}}\right)^{1/2}$
                
                \item Huber function
                \begin{equation}
                    h_{\delta}(x) = 
                    \begin{cases}
                        \frac{L}{2}x^2 & \text{if } x \leq \delta \\
                        L \delta x - \frac{L}{2}\delta^2 & \text{if } x > \delta
                    \end{cases}
                    \label{eq:huber_functions}
                \end{equation}
        
                \item For $i \in [|0, n-1|]$, define $\xi_i \triangleq \delta \left( 1 + \sum_{k=i}^{n-2} L \gamma_i \right)$.
                
                \item $\lambda = \frac{L\eta}{(1 + L\gamma_{n-2})(1 + \eta^2 + \xi_0^2)}$.
            \end{itemize}
            
            Based on those notations, we define the $3$-dimensional function
            \begin{equation}
                f(x) = \max\left[\frac{L}{2}\left(x^{(1)} - 1 + |x^{(2)}|\sqrt{3}\right), h_{\delta}\left(x^{(3)}\right), \frac{\lambda}{2}\|x\|_2^2\right].
            \end{equation}
            
            $f$ is convex as maximum of 3 convex functions.
            
            Moreover, we note that $\mathcal{X}_\star = \lbrace 0 \rbrace$ and $f_\star = 0$.
            
            And each of the three components defining $f$ is smaller that $\frac{L}{2}\|x\|_2^2$.
            Indeed,
            \begin{align}
                \frac{L}{2}\left(x^{(1)} - 1 + |x^{(2)}|\sqrt{3}\right) & = \frac{L}{2}\left(\left(x^{(1)}\right)^2 - \left(x^{(1)} - \frac{1}{2}\right)^2 + \left(x^{(2)}\right)^2 - \left(x^{(2)} - \frac{\sqrt{3}}{2}\right)^2 \right) \nonumber \\
                & \leq \frac{L}{2} \left( \left(x^{(1)}\right)^2 + \left(x^{(2)}\right)^2 \right) \nonumber \\
                & \leq \frac{L}{2} \|x\|_2^2 \\
                h_{\delta}\left(x^{(3)}\right) & \leq \frac{L}{2}(x^{(3)})^2 \nonumber \\
                & \leq \frac{L}{2}\|x\|_2^2 \\
                \lambda & \leq \frac{L\eta}{2\eta} = \frac{L}{2} \leq L, \text{ hence }
                \frac{\lambda}{2}\|x\|_2^2 \leq \frac{L}{2}\|x\|_2^2
            \end{align}
            
            Therefore, $f$ is also $\QG^+(L)$.
            
            We choose to start the GD algorithm at $x_0 \triangleq \begin{pmatrix} 1 & \eta & \xi_0 \end{pmatrix}^\top$.
            
            We claim that after $i$ ($ 0 \leq i \leq n-1$) steps of GD, $x_i = \begin{pmatrix} 1 & \eta & \xi_i \end{pmatrix}^\top $.
            
            This can be proven by induction.
            Indeed, by definition, this is true for $i=0$.
            We now assume this property is true for some $i < n-1$ and want to prove it for $i+1$.
            
            From the 3 remarks
            \begin{align}
                h_{\delta}(\xi_i) & \geq \frac{L}{2}\left(1 - 1 + |\eta|\sqrt{3}\right) \\
                h_{\delta}(\xi_i) & \geq \frac{\lambda}{2}\|x_i\|_2^2 \\
                \xi_i & \geq \delta,
            \end{align}
            
            we conclude that $\nabla f(x_i) = \begin{pmatrix} 0 & 0 & L\delta \end{pmatrix}^\top$.
            
            Hence $x_{i+1} = x_i - \begin{pmatrix} 0 & 0 & L \gamma_i \delta \end{pmatrix}^\top = \begin{pmatrix} 1 & \eta & \xi_{i+1} \end{pmatrix}^\top$.
            
            Finally, from the 2 remarks
            \begin{align}
                \frac{L}{2}(1 - 1 + |\eta|\sqrt{3}) & \geq h_{\delta}(\xi_{n-1}) \\
                \frac{L}{2}(1 - 1 + |\eta|\sqrt{3}) & \geq \frac{\lambda}{2}\|x_{n-1}\|_2^2,
            \end{align}
            
            we conclude that $\nabla f(x_{n-1}) = \frac{L}{2} \begin{pmatrix} 1 & \sqrt{3} & 0 \end{pmatrix}^\top$, leading to $x_n = x_{n-1} - \gamma_{n-1} \frac{L}{2} \begin{pmatrix} 1 & \sqrt{3} & 0 \end{pmatrix}^\top = \begin{pmatrix} 1 - \frac{L\gamma_{n-1}}{2} & \eta - \frac{L\gamma_{n-1}\sqrt{3}}{2} & \delta \end{pmatrix}^\top$.
            
            We compute the two quantities
            \begin{align}
                \|x_0\|^2 & = 1 + \eta^2 + \xi_0^2 \\
                f(x_n) & \geq \frac{L}{2}\left(1 - \frac{L\gamma_{n-1}}{2} - 1 + |\eta - \frac{L\gamma_{n-1}\sqrt{3}}{2}|\sqrt{3}\right) \nonumber \\
                & = \frac{L}{2}\left(- \frac{L\gamma_{n-1}}{2} + \left(\frac{L\gamma_{n-1}\sqrt{3}}{2} - \eta\right)\sqrt{3}\right) \nonumber \\
                & = \frac{L}{2}\left(L\gamma_{n-1} - \eta\sqrt{3}\right)
            \end{align}
            
            Finally,
            
            \begin{align}
                \frac{f(x_n) - f_\star}{d(x_0, \mathcal{X}_\star)^2} & \geq \frac{L}{2} \frac{L\gamma_{n-1} - \eta\sqrt{3}}{1 + \eta^2 + \xi_0^2} \nonumber \\
                & = \frac{L}{2} \frac{L\gamma_{n-1} - \eta\sqrt{3}}{1 + \eta^2 + \delta^2 \left( 1 + \sum_{k=i}^{n-2} L \gamma_i \right)^2} \nonumber \\
                & = \frac{L}{2} \frac{L\gamma_{n-1} - \eta\sqrt{3}}{1 + \eta^2 + \frac{\eta\sqrt{3}}{1+L\gamma_{n-2}} \left( 1 + \sum_{k=i}^{n-2} L \gamma_i \right)^2} \nonumber \\
                & \underset{\eta \rightarrow 0}{\longrightarrow} \frac{L}{2}L\gamma_{n-1}.
            \end{align}
            
            Hence, for any $\epsilon >0$, we can find $\eta >0$ sufficiently small such that $f$ reaches the claim of the Theorem.
            
        $\hfill\blacksquare$

    \subsection{A new tuning prescription}
    \label{apx:gd_conjecture}

        Theorem~\ref{thm:non_convergence_gd} provides a lower bound on the last iterate value of the subgradient method on the class of $\QG^+$ convex functions.
        Moreover, a new analysis of the subgradient method on the Huber function~\eqref{eq:huber_functions}, starting at $x_0 = 1 + 2 \sum_{k=0}^{n-1} L\gamma_k$ provides another lower bound.
        
        Combining those 2 results, we know that whatever $(\gamma_i)_{0 \leq i \leq n-1}>0$ is, there exists $f$ an $L-\QG^+$ convex function as well as a starting point $x_0$ such that
        
        \begin{equation}
            \label{eq:gd_lower_bound}
            f(x_n) - f_\star \geq \max\left(\frac{L}{2}\frac{1}{1 + 2 \sum_{k=0}^{n-1} L\gamma_k}, \frac{L}{2}L\gamma_{n-1}\right) d(x_0, \mathcal{X}_\star)^2.
        \end{equation}

        Naturally, we propose the sequence of
        
        \begin{wrapfigure}[9]{R}{0.55\textwidth}
            \vspace{-1cm}
            \begin{minipage}{0.45\textwidth}
                \begin{algorithm}[H]
                    \caption{GD with decreasing step-sizes}\label{alg:gd_decreasing_step_sizes}
                    \KwInput{$x_0$, $L$}
                    $u_0 = 1$
                    
                    \For{k=1 \ldots n}{
                        $u_{k} \gets \frac{u_{k-1}}{2} + \sqrt{\left(\frac{u_{k-1}}{2}\right)^2 + 2}$;
                        
                        $\gamma_{k-1} \gets \frac{1}{L u_k}$;
                        
                        \text{Pick } $g_{k-1} \in \partial f(x_{k-1})$;
                        
                        $x_k \gets x_{k-1} - \gamma_{k-1} g_{k-1}$
                    }
                    \KwOutput{$x_n$}
                \end{algorithm}
            \end{minipage}
        \end{wrapfigure}

        \noindent $(\gamma_i)_{0 \leq i \leq n-1}>0$ that verifies for all $n$, $\frac{L}{2}\frac{1}{1 + 2 \sum_{k=0}^{n-1} L\gamma_k} = \frac{L}{2}L\gamma_{n-1}$ (for each index $n \geq 1$).
        This is summarized in Algorithm~\ref{alg:gd_decreasing_step_sizes}.
        We note that $u_n \sim 2\sqrt{n}$ and $\gamma_{n-1} \sim \frac{1}{2L\sqrt{n}}$.
        The lower bound~\eqref{eq:gd_lower_bound} for this method becomes $f(x_n) - f_\star \geq \frac{L}{2}L\gamma_{n-1} d(x_0, \mathcal{X}_\star)^2 \sim \frac{L}{4\sqrt{n}} d(x_0, \mathcal{X}_\star)^2$.
        We conjecture that this bound is actually reached by the proposed method~\ref{alg:gd_decreasing_step_sizes}.

        \begin{restatable}{Conj}{convergence_of_gd_with_decreasing_step-sizes}\textbf{\emph{(Convergence of GD with decreasing step-sizes)}}
            \label{conj:gd_sqrt}
            The algorithm~\ref{alg:gd_decreasing_step_sizes} verifies the following lower bound on every $L-\QG^+$ convex function $f$:
            \begin{equation}
                f(x_n) - f_\star \leq \frac{L}{2 u_n} d(x_0, \mathcal{X}_\star)^2.
            \end{equation}
            where $u_n$ is the sequence used in~\ref{alg:gd_decreasing_step_sizes}, defined by
            \begin{align}
                u_0 & = 1 \\
                u_{k} & = \frac{u_{k-1}}{2} + \sqrt{\left(\frac{u_{k-1}}{2}\right)^2 + 2}, \quad \text{for every } k \in \llbracket 1, n \rrbracket.
            \end{align}
            and verifying
            \begin{equation}
                u_n \sim 2\sqrt{n}.
            \end{equation}
        \end{restatable}

\begin{figure}[h!]
            \centering
            \includegraphics[width=0.9\textwidth]{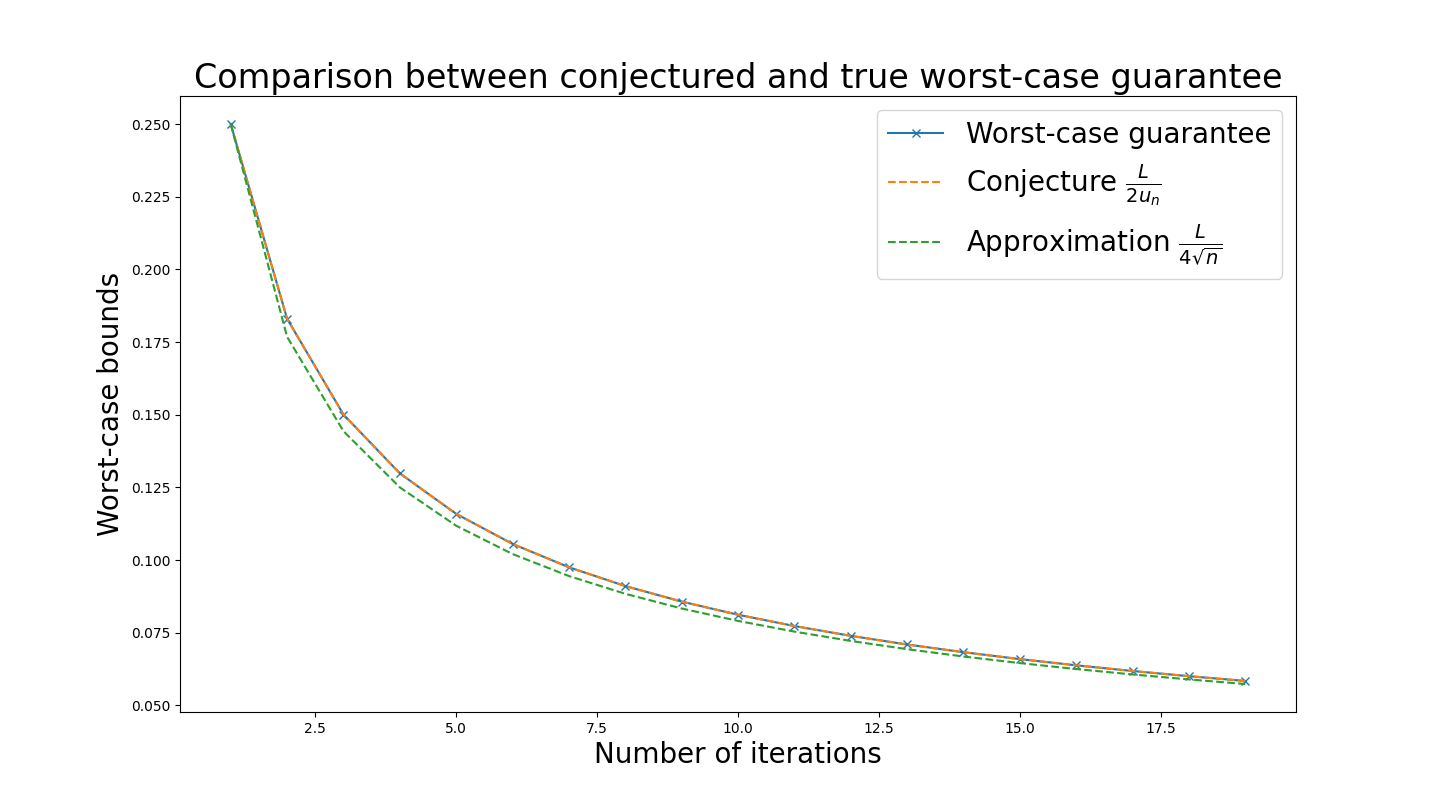}
            \caption{Verification of the conjecture using PEPs}
            \label{fig:gd_sqrt}
        \end{figure}

        This conjecture is supported by the Figure~\ref{fig:gd_sqrt} that has been built using the PEP framework.
        This figure represents the worst-case guarantee of Algorithm~\ref{alg:gd_decreasing_step_sizes} as a function of the number of iterations.
        The conjecture (red curve) follows exactly the numerical worst-case guarantee provided by the PEPs (blue curve) and the equivalent sequence (green curve) is very close to the 2 previous ones.
        
\section{First-order lower bound}
\label{apx:lower_bound}

    In this section, we prove the lower bound of first order methods on $\QG^+$ convex functions.
    
    In the first subsection, we assume that the iterates of the first order algorithm must stay in the span of the past observed gradients.
    
    In the second subsection, we release this assumption and still prove the same lower bound.
    This proof is a bit more technical, hence the reason why we provide the two proofs.
    
    \subsection{\texorpdfstring{Proof of Theorem 2.3}{Proof of Theorem~\ref{thm:general_lower_bound}}}
    \label{apx:lower_bound_1}

        The following theorem brings a lower bound over all first order methods verifying that all the iterates lie into the span of the previously observed gradients.
        
        \lowerboundoffirstorderalgorithm*
        
        \noindent \textit{Proof.}
            Consider $f(x) \triangleq \frac{L}{2}\|x\|_\infty^2$ defined on $\mathbb{R}^{n+1}$ and $x_0=\Vec{\mathbf{1}}$.
            
            After $k$ steps, the oracle can ``choose'' to return a vector that lies in the first $k+1$ dimension of the input space, leading to $f(x_n) - f_\star \geq f(x_0) - f_\star = \frac{L}{2} = \frac{L}{2}\frac{1}{n+1} d(x_0, \mathcal{X}_\star)^2$.
        $\hfill\blacksquare$

        \begin{Rem}
            Note that by considering instead $x_0 = \begin{pmatrix}
                1 & 1 - \varepsilon/n & 1 - 2\varepsilon/n & \ldots & 1 - \varepsilon
            \end{pmatrix}^\top$,
            one ends up with $f(x_n) - f_\star \geq \frac{L}{2}\frac{1 - \varepsilon}{n+1} d(x_0, \mathcal{X}_\star)^2$ \textbf{whatever} the oracle ``chooses'' to return.
        \end{Rem}

    \subsection{Lower bound proof without span assumption}
    \label{apx:lower_bound_2}

        In this section we release the span assumption and prove that the previously shown lower bound still holds.
        This proof is a bit more technical than the one of Theorem~\ref{thm:general_lower_bound} proven in Appendix~\ref{apx:lower_bound_1}.
        
        \begin{Th}[Lower bound of first order algorithm without span assumption]
            Let $\mathcal{A}$ a first-order algorithm.
            Then, for any $n \geq 0$, there exists $d$ a positive integer, $f$ a $L-\QG^+$ convex function of input space $\mathbb{R}^d$ and a starting point $x_0$ such that $f(x_n) - f_\star \geq \frac{L}{2}\frac{1}{n+1} d(x_0, \mathcal{X}_\star)^2$.
        \end{Th}
        
        \noindent \textit{Proof.}
            Let $\mathcal{E} = \left\{ v \in \mathbb{R}^{n+1} | \forall i \in \llbracket 1, n+1 \rrbracket, |v_i| = 1 \right\}$ the set of the $2^{n+1}$ vectors of $\mathbb{R}^{n+1}$ which all coordinates are $\pm 1$.
            
            For each $v\in\mathcal{E}$, we introduce $f_v(x) \triangleq \frac{L}{2}\|x - v\|_\infty^2$ defined on $\mathbb{R}^{n+1}$.
            
            First note that all those functions are $L-\QG^+$ convex.
            We will prove that not only there exists a starting point $x_0$ and a $L-\QG^+$ convex function such that $f(x_n) - f_\star \geq \frac{L}{2}\frac{1}{n+1} d(x_0, \mathcal{X}_\star)^2$, but also that there exists a starting point $x_0$ and a function among the $f_v$ we introduced above such that the latest holds.
            
            To proceed, we need to show that the algorithm $\mathcal{A}$ cannot know the right $v$ after only $n$ iterations and therefore, cannot guarantee $f(x_n) - f_\star \leq \frac{L}{2}$.
            Taking $x_0 = \Vec{0}$, $\|x_0 - x_\star \|^2 = n+1$ whatever $x_\star$ is (among $\mathcal{E}$), hence the result.
            
            In order to prove that the algorithm cannot know the solution after $n$ iterations, we keep track of all the remaining possibilities across time.
            
            We denote by $\mathcal{E}_k$ the remaining possibilities after $k$ steps of the algorithm.
            In particular, $\mathcal{E}_0=\mathcal{E}$.
            
            At each step $k$, $\mathcal{A}$ guesses $x_k$ based on all the previous information, summarized in $\mathcal{E}_k$, and the oracle provides $f(x_k)$ and a sub-gradient $g_k \in \partial f(x_k)$.
            
            Let $v_k \in \arg\max_{v \in \mathcal{E}_k} \|x_k - v_k\|_\infty$.
            We consider the case where the oracle behaves like if the objective function to minimize was $f_{v_k}$.
            Moreover, in the case where $f_{v_k}$ is not differentiable in $x_k$, we ask that the oracle returns a sub-gradient co-linear to a vector belonging to the canonical basis (which is always possible).
            
            Considering $i$ such that $g_k$ is co-linear to $e_i$, we obtain
            $\mathcal{E}_{k+1} = \left\{ v \in \mathcal{E}_k | \left< v, e_i \right> = \left< v_k, e_i \right> \right\}$, reducing by half the number of remaining elements at each step (except when the algorithm badly guesses and receive twice the same direction, in which case one of the steps is useless).
            
            After $n$ steps, $\mathcal{E}_{n}$ contains 2 elements, and the algorithm $\mathcal{A}$ must guess based on nothing.
            Again, we consider the further one to the last guess as the right solution, and obtain the lower bound provided by the Theorem.
            
        $\hfill\blacksquare$

\section{Main result: worst-case guarantee of proposed methods}
\label{apx:main_result}

    In this section, we prove Theorem~\ref{thm:main}, the main result of this paper, stating that all the sequences of iterates verifying a certain property enjoy an upper bound guarantee corresponding to the lower bound presented in \Cref{thm:general_lower_bound} and proved in \Cref{apx:lower_bound}.
    
    \mainresult*
    
    \noindent \textit{Proof.}
        This proof relies on the Lyapunov function
        \begin{equation}
            V_n \triangleq n(f(x_{n-1}) - f_\star) + \frac{L}{2} \left\| x_0 - \pi_{\mathcal{X_\star}}(x_0) - \sum_{i=0}^{n-1}\frac{1}{L}g_i \right\|^2.
        \end{equation}
        
        For all $k$, we verify
        
        \begin{eqnarray*}
            V_{k+1} - V_k & = & \left[(k+1)(f(x_{k}) - f_\star) + \frac{L}{2} \left\| x_0 - \pi_{\mathcal{X_\star}}(x_0) - \sum_{i=0}^{k}\frac{1}{L}g_i \right\|^2 \right] \\
            & & - \left[k(f(x_{k-1}) - f_\star) + \frac{L}{2} \left\| x_0 - \pi_{\mathcal{X_\star}}(x_0) - \sum_{i=0}^{k-1}\frac{1}{L} g_i \right\|^2 \right] \\
            & = & (f(x_k) - f_\star) + k(f(x_k) - f(x_{k-1})) + \frac{L}{2} \left[ -\frac{2}{L} \left< g_k , x_0 - \pi_{\mathcal{X_\star}}(x_0) - \sum_{i=0}^{k-1}\frac{1}{L} g_i \right> + \frac{1}{L^2} \|g_k\|^2 \right] \\
            & = & \left(f(x_k) - f_\star + \frac{1}{2L} \|g_k\|^2\right) + k(f(x_k) - f(x_{k-1})) - \left< g_k , x_0 - \pi_{\mathcal{X_\star}}(x_0) - \sum_{i=0}^{k-1}\frac{1}{L} g_i \right> \\
            & \overset{(\ref{eq:interp_convexity},~\ref{eq:interp_convexity_qg})}{\leq} & \left< g_k , x_k - \pi_{\mathcal{X_\star}}(x_0) \right> + k\left< g_k , x_k - x_{k-1} \right> - \left< g_k , x_0 - \pi_{\mathcal{X_\star}}(x_0) - \sum_{i=0}^{k-1}\frac{1}{L} g_i \right> \\
            & = & \left< g_k , x_k - \pi_{\mathcal{X_\star}}(x_0) + k(x_k - x_{k-1}) - \left( x_0 - \pi_{\mathcal{X_\star}}(x_0) - \sum_{i=0}^{k-1}\frac{1}{L} g_i \right) \right> \\
            & = & \left< g_k , (k+1) x_k - k x_{k-1} - x_0 + \sum_{i=0}^{k-1}\frac{1}{L} g_i \right>
        \end{eqnarray*}
        
        The assumption therefore concludes
        \begin{equation*}
            \forall k, V_{k+1} \leq V_k
        \end{equation*}
        
        Finally,
        
        \begin{equation*}
            (N+1)(f(x_{N}) - f_\star) + \frac{L}{2} \left\| x_0 - \pi_{\mathcal{X_\star}}(x_0) - \sum_{i=0}^{N}\frac{1}{L} g_i \right\|^2 = V_{n+1} \leq V_0 = \frac{L}{2}\|x_0 - x^*\|^2.
        \end{equation*}
        
        Hence,
        
        \begin{equation*}
            f(x_{N}) - f_\star \leq \frac{L}{2}\frac{1}{N+1}\|x_0 - x^*\|^2.
        \end{equation*}
        
    $\hfill\blacksquare$

\section{\texorpdfstring{Summary of convergence results on $\QG^+$ convex and Lipschitz convex}{Summary of convergence results on QG+ convex and Lipschitz convex}}
\label{apx:tab}

    In this section, we state and prove the 2 results of \Cref{tab:optimality_summary} that are not already proven elsewhere.
    
    \addtocounter{table}{-1}
    
    \begin{table*}
        {
        \caption{
            Optimality of the proposed methods over the set of $\QG^+$ convex functions and $M$-Lipschitz convex functions.
            ELS: Exact Line-Search. $\cmark$ indicates optimality among the class and $\xmark$ the contrary.
            All counter examples are given in App.~\ref{apx:tab}. $\ ^{\dagger}$:~constants resulting in optimal convergence rates depend on the class, thus for example heavy-ball with step-size $\frac{\text{constant}}{(t+2)}$ is not adaptive as it does not achieve the optimal rate for both classes with the same constant. $\ ^\ddagger$: up to a $\log$ factor.
            } \vspace{-0.3cm}
        \begin{center}
            {\renewcommand{\arraystretch}{1.8}
            \resizebox{\linewidth}{!}{
             \begin{tabular}{@{}lllcrcrc@{}}
                \specialrule{2pt}{1pt}{1pt}
                \multicolumn{3}{c}{Method} & \multicolumn{4}{c}{Function class} & Parameter free \\
                 \cmidrule{1-3}\cmidrule(l){4-7}\cmidrule(l){8-8}
                Algorithm & Step-sizes $(\gamma_t)_{0\le t\le n-1}$  & Iterate & \multicolumn{2}{c}{$\QG^+(L)$ convex} & \multicolumn{2}{c}{$M$-Lipschitz convex}   \\
                \cmidrule(l){4-5}\cmidrule(l){6-7}
                Subgradient (Alg.~\ref{alg:subgrad}) & $\text{constant}^{\dagger}$ & Average &  $\cmark$ & (Thm.~\ref{thm:gd_average}) & $ \xmark$ & (Thm.~\ref{thm:subgrad_constant_lower_lip}) & $\xmark$ \\
                Subgradient (Alg.~\ref{alg:gd_decreasing_step_sizes}) & $ \text{constant}^{\dagger}/\sqrt{t}$ & Average & $\xmark$ &~\eqref{eq:LB_smooth} & $\cmark^\ddagger$ &~\citep[Sec. 3.2.3]{Nest03a} & $\xmark$ \\
                Subgradient (Alg.~\ref{alg:subgrad_els}) & ELS & Average &  $\xmark$ & (Thm.~\ref{thm:gd_els_lower_bound}) & $\xmark$ & (Thm.~\ref{thm:gd_els_lower_bound}) & $\cmark$ \\
                Subgradient (Alg.~\ref{alg:subgrad_els}) & ELS & Last &  $\xmark$ & (Thm.~\ref{thm:gd_els_lower_bound}) & $\xmark$ & (Thm.~\ref{thm:gd_els_lower_bound}) & $\cmark$ \\
                \midrule
                Heavy-ball (Alg.~\ref{alg:ogm}) & $\text{constant}^{\dagger}/(t+2)$ & Last & $\cmark$ & (Cor.~\ref{cor:optimal}) & $\cmark$ &~\citep[][, Cor. 3]{drori2020efficient} & $\xmark$ \\
                Heavy-ball (Alg.~\ref{alg:ogm_ls}) & ELS & Last & $\cmark$ & (Cor.~\ref{cor:optimal}) & $\cmark$ &~\citep[][, Cor. 4]{drori2020efficient} & $\cmark$ \\
                \specialrule{2pt}{1pt}{1pt}\vspace{0em}
            \end{tabular}}
            \vspace{-1cm}
            }
        \end{center}}
    \end{table*}
    
    We first state \Cref{thm:subgrad_constant_lower_lip}.
    
    \begin{Th}
        For any $M>0$ and any $\gamma>0$, the subgradient method~\ref{alg:subgrad} with constant step-size $\gamma$ cannot be guaranteed to converge to optimum on all the $M$-Lipschitz continuous convex functions, both in last iterate and in Polyak-Rupert averaged iterate.
        \label{thm:subgrad_constant_lower_lip}
    \end{Th}

    \noindent \textit{Proof.}
        First we note that $f \triangleq z \mapsto M |z|$ is $M$-Lipschitz continuous and convex.
        Let $\gamma>0$.
        We consider $x_0 = \frac{3}{4}M\gamma$ the starting point of the subgradient method with constant step-size $\gamma$.
        We verify that $x_1 = -\frac{1}{4}M\gamma$ and that the sequence $(x_t)_t$ cycles back to $\frac{3}{4}M\gamma$.
        Therefore, the sequence itself does not converge and the sequence of the PR averaged iterates converges to $\frac{1}{4}M\gamma$, while the optimum value would be 0.
    $\hfill\blacksquare$

    \begin{wrapfigure}[6]{R}{0.56\textwidth}
        \begin{minipage}{0.56\textwidth}
            \vspace{-1.0cm}
            \begin{algorithm}[H]
                \caption{Subgradient method with line-search \label{alg:subgrad_els}}
                \KwInput{$x_0$, $v_0 \gets 0$}
                \For{$k=1 \ldots n$}{
                    \text{Pick } $g_{k-1} \in \partial f(x_{k-1})$.
                    
                    $\alpha_k \gets \arg\min_{\alpha} f\left( x_{k-1} - \alpha g_{k-1} \right)$
                    
                    $x_k \gets x_{k-1} - \alpha_k g_{k-1}$
                }
                \KwOutput{$x_n$}
            \end{algorithm}
        \end{minipage}
    \end{wrapfigure}

    Then, one could wonder whether performing exact line search steps on the subgradient method (Algorithm~\ref{alg:subgrad_els}.) leads to convergence on Lipschitz continuous convex or $\QG^+$ convex function.
    We now prove Theorem~\ref{thm:gd_els_lower_bound} stating that the subgradient method with exact line search does not converge for all functions neither of the class of Lipschitz continuous convex functions nor of the class of $\QG^+$ convex functions, neither in last iterate nor in Polyak-Rupert averaged iterate.
    
    \begin{Th}
        There exists a $\QG^+$ convex function $f_{QG}$ and a Lipschitz continuous convex function $f_{Lip}$ such that, the iterates $((x_n)_n^{f_{QG}})_n$ and $((x_n)_n^{f_{Lip}})_n$ obtained by Algorithm~\ref{alg:subgrad_els}, verify the 2 guarantees
        
        \begin{align}
            f(x_n) - f_\star & \geq && \frac{L}{6} \|x_0 - x_\star\|^2. \\
            f(\bar{x_n}) - f_\star & \geq && \frac{L}{6} \|x_0 - x_\star\|^2.
        \end{align}
        
        where $\bar{x_n}$ denotes the Polyak-Rupert averaged iterate obtained from the sequence $(x_n)_n$.
        
        \label{thm:gd_els_lower_bound}
    \end{Th}

    \noindent \textit{Proof.}
        Considering $f_{Lip}(z) = M\|x\|_\infty$ and $f_{QG}(z) = \frac{L}{2}\|x\|_\infty^2$ defined on $\mathbb{R}^3$,
        the iterates of Algorithm~\ref{alg:subgrad_els} can be cycling between the four points $(1,1,1)$, $(1,-1,1)$, $(-1,1,1)$ and $(-1,-1,1)$.
        Therefore, the aforementioned statement.
    $\hfill\blacksquare$
    
\vspace{.5cm}

\section{\texorpdfstring{Interpolation results for $\QG^+$ convex functions}{Interpolation results for QG+ convex functions}}
\label{apx:interpolation_conditions}

    The interpolation conditions of a given class represent the key ingredient to use the PEP framework on this class.
    Theorem~\ref{thm:interp} provides the interpolation conditions for the class of $\QG^+$ convex functions.
    In this section, we recall this result and prove it.

    \interpolationconditions*

    \noindent \textit{Proof.}

        We prove the two implications one by one.

        \begin{itemize}
            \item[$\Rightarrow$:] Assume there exists such a convex-$\QG^+$ function $f$ that interpolates $(x_i, g_i, f_i)_{i \in I}$.
            Equation~\eqref{eq:interp_convexity} follows immediately from convexity.
            Let's prove equation~\eqref{eq:interp_convexity_qg}.
            Let $i \in I_\star, \forall j \in I$ and $x \in \mathbb{R}^d$.
            We have:
            \begin{equation*}
                f_j + \left< g_j, x - x_j \right> \overset{\text{CVX}}{\leq} f(x) \overset{\QG^+}{\leq} \min_{z \in \mathbb{R}^d} f(z) + \frac{L}{2}d(x, \mathcal{X}_\star)^2 \leq f_i + \frac{L}{2}\|x - x_i \|^2.
            \end{equation*}
            Rewriting the previous equation for $x = x_i + \frac{1}{L}g_j$ leads to equation~\eqref{eq:interp_convexity_qg}.

            \item[$\Leftarrow$:] Let's consider equations~\eqref{eq:interp_convexity} and~\eqref{eq:interp_convexity_qg} are verified.
            Applying~\eqref{eq:interp_convexity} with $j \in I_\star$
            \begin{equation}
                \forall i \in I, \forall j \in I_\star, f_i \geq f_j \label{eq:stationary_are_minima}
            \end{equation}
            In particular, $\forall i \in I_\star, \forall j \in I_\star, f_i = f_j$.
            Hence, let's introduce $f_\star$ the common value of all the $f_i$ for $i \in I_\star$.
            Let's denote here $\mathcal{X}^\star$ the convex hull of $\left\{x_i\right\}_{i \in I_\star}$.
            Finally let's introduce $\mu \triangleq 2 \min_{i \in I \setminus I_\star} \left( \frac{f_i - f_\star}{d(x_i, \mathcal{X}^\star)^2} \right)$.

            Let's prove that the following function $f$ is a solution:
            \begin{equation}
                f(x) = \max\left( \max_{j \in I} \left( f_j + \left< g_j, x - x_j \right> \right), f_\star + \frac{\mu}{2}d\left(x, \mathcal{X}^\star \right)^2 \right).
                \label{eq:interpolation_function}
            \end{equation}

            \begin{itemize}
                \item[-]\underline{$\forall i \in I, f(x_i) = f_i$:} For all $i \in I$, equation~\eqref{eq:interp_convexity} shows $\max_{j \in I} \left( f_j + \left< g_j, x_i - x_j \right> \right) \leq f_i$ and the definition of $\mu$ leads to $f_\star + \frac{\mu}{2}d\left(x_i, \mathcal{X}^\star \right)^2 \leq f_i$.
                Hence, for all $i \in I$, $f(x_i) \leq f_i$.
                Moreover, from equation~\eqref{eq:interpolation_function}, $f(x) \geq \left( f_i + \left< g_i, x - x_i \right> \right)$, hence $f(x_i) \geq f_i$.
                Finally, we conclude $\forall i \in I, f(x_i) = f_i$.

                \item[-]\underline{$\forall i \in I, g_i \in \partial f(x_i)$:} $\forall i \in I, f(x) \geq \left( f_i + \left< g_i, x - x_i \right> \right)$, and from the previous point, we conclude $\forall i \in I, f(x) \geq \left( f(x_i) + \left< g_i, x - x_i \right> \right)$.
                Finally, $\forall i \in I, g_i \in \partial f(x_i)$.

                \item[-]\underline{$f$ is convex:} $f$ is clearly defined as the maximum of convex functions, hence is convex.

                \item[-]\underline{$f$ is $\QG^+$:} We aim at proving that $\forall x \in \mathbb{R}^d, f(x) \leq f_\star + \frac{L}{2}d\left(x, \mathcal{X}^\star \right)^2$.
                Since it is clear that $\forall x \in \mathbb{R}^d, f_\star + \frac{\mu}{2}d\left(x, \mathcal{X}^\star \right)^2 \leq f_\star + \frac{L}{2}d\left(x, \mathcal{X}^\star \right)^2$, it remains to prove that
                $\forall x \in \mathbb{R}^d, \forall j \in I, f_j + \left< g_j, x - x_j \right> \leq f_\star + \frac{L}{2}d\left(x, \mathcal{X}^\star \right)^2 $.
                The latest is also equivalent to $\forall x \in \mathbb{R}^d, \forall j \in I, \forall x_\star \in \mathcal{X}^\star, f_j + \left< g_j, x - x_j \right> \leq f_\star + \frac{L}{2}\| x - x_\star \|^2 $.
                For $j$ and $x_\star$ fixed, this expression is a quadratic form in $x$, optimized for $x = x_\star + \frac{1}{L}g_j$.
                Hence, we need to show $\forall j \in I, \forall x_\star \in \mathcal{X}^\star, f_j + \left< g_j, x_\star + \frac{1}{L}g_j - x_j \right> \leq f_\star + \frac{L}{2}\left\| x_\star + \frac{1}{L}g_j - x_\star \right\|^2 $, also rewritten $\forall x_\star \in \mathcal{X}^\star, \forall j \in I, f_\star \geq f_j + \left< g_j, x_\star - x_j \right> + \frac{1}{2L} \|g_j\|^2$.
                Since $\mathcal{X}^\star$ is the convex hull of $\left\{ x_i \right\}_{i \in I_\star}$, the latest is obtained by linear combination (with non-negative weights) of equation~\ref{eq:interp_convexity_qg} for different values of $i$.
            \end{itemize}
        \end{itemize}
    $\hfill\blacksquare$

\section{Convergence bound on other classes}
\label{apx:upper_assumption}

    In this section, we naturally extend the previous results to the class of $h-\RG^+$ (See. \Cref{def:h_rg}) convex functions.

    \hbgeneral*

    \noindent \textit{Proof.}

        First note that since $h$ is strictly increasing, $h^{-1}$ is well defined.
        Furthermore, Definition~\ref{def:h_rg} can be expressed as
        $h^{-1}\left(f(x)-f_\star\right) \leq d(x, \mathcal{X}_\star)^2.$
        Therefore, $h^{-1}\left(f-f_\star\right)$ is $2-\QG^+$.
        And we know by definition that $h$ is increasing and concave, then $h^{-1}$ is increasing and convex.
        Since $f$ is also convex, so is $h^{-1}\left(f-f_\star\right)$.

        We conclude that $h^{-1}\left(f-f_\star\right)$ is $2-\QG^+$ and convex, and then we know that Algorithm~\ref{alg:hb_ogm} applied on $h^{-1}\left(f-f_\star\right)$ leads to

        \begin{equation*}
            h^{-1}\left(f(x_n) - f_\star\right) \leq \frac{d(x_0, \mathcal{X}_\star)^2}{n+1}.
        \end{equation*}

        It remains to compose the above by $h$ and to notice that Algorithm~\ref{alg:hb_ogm} applied on $h^{-1}\left(f-f_\star\right)$ is exactly Algorithm~\ref{alg:hb_general}.

    $\hfill\blacksquare$

\section{Linear convergence guarantees under lower bound assumption}
\label{apx:restart}

    In all this section, we assume that $f$ is convex and $h-\RG^+$ for a certain $h$.
    Moreover, we consider that $f$ verifies the following additional assumption (referred to as ``\L{}ojasiewicz error bound inequality'' in~\citep{bolte2017error}) for a given $\kappa \geq 1$:

    \begin{Assump}
        \label{assump:h_kappa}
        For all $x \in \mathbb{R}^d$, $f(x)-f_\star \geq h\left( \frac{ d(x, \mathcal{X}_\star)^2}{\kappa} \right)$.
    \end{Assump}

    \begin{Rem}
        When $h$ is the linear function $h: R \mapsto \frac{LR}{2}$, then $f$ is simply $L-\QG^+$ convex as well as $\mu-\QG^-$ $\left(\text{where } \mu \triangleq \frac{L}{\kappa}\right)$(See~\citep{guille2021study} for the definition of $\QG^-$).
    \end{Rem}

    We introduce the following algorithm based on the restart idea studied in~\citep{nemirovskii1985optimal, nesterov2013gradient, iouditski2014primal}.

    \begin{algorithm}
        \caption{Heavy-ball with restart}
        \label{alg:hb_restart}
        \KwInput{$x_0$, $h$, $f_\star$}
        \For{$k=1 \ldots n$}{
        Choose ~ $g_{k-1}$ from $\partial f(x_{k-1})$

        $l \gets k \text{ mod } \left\lfloor \kappa e \right\rfloor - 1 $ \quad (between $1$ and $\left\lfloor \kappa e \right\rfloor - 1$).

        $x_k \gets x_{k-1} - \frac{1}{2(l+1)}\frac{1}{h' \circ h^{-1}\left(f(x_{k-1}) - f_\star\right)} g_{k-1} + \frac{l-1}{l+1}\left( x_{k-1} - x_{k-2} \right) $
        }
        \KwOutput{$x_n$}
    \end{algorithm}

    This algorithm comes with the linear convergence rate guarantee

    \begin{center}
        \fbox{\parbox{\textwidth}{
        \begin{restatable}{Th}{hbrestart}\textbf{\emph{}}
            \label{thm:hb_restart}
            Algorithm~\ref{alg:hb_restart} verifies for every $n$ multiple of $\left\lfloor \kappa e \right\rfloor - 1$:

            \begin{equation}
                d(x_n, \mathcal{X}_\star)^2 \leq \left(1 - \frac{1}{\kappa e}\right)^n d(x_0, \mathcal{X}_\star)^2.
            \end{equation}
        \end{restatable}
        }}
    \end{center}

    \noindent \textit{Proof.}

        From Theorem~\ref{thm:hb_general}, running Algorithm~\ref{alg:hb_general} leads to the guarantee
        \begin{equation}
            f(x_n) - f_\star \leq h\left(\frac{d\left(x_0, \mathcal{X}_\star\right)^2}{n+1} \right).
        \end{equation}

        From the additional assumption (\ref{assump:h_kappa}), we can upper bound the left hand size of the above and write
        \begin{equation*}
            h\left( \frac{ d(x_n, \mathcal{X}_\star)^2}{\kappa} \right) \leq h\left(\frac{d\left(x_0, \mathcal{X}_\star\right)^2}{n+1} \right).
        \end{equation*}

        Hence,
        \begin{equation*}
            d(x_n, \mathcal{X}_\star)^2 \leq \frac{\kappa}{n+1}d(x_0, \mathcal{X}_\star)^2.
        \end{equation*}

        The latest is a contraction guarantee.
        Indeed, for $n$ sufficiently large, $d(x_n, \mathcal{X}_\star)^2 < d(x_0, \mathcal{X}_\star)^2$.
        The average contraction factor is $\left(\frac{\kappa}{n+1}\right)^{1/n}$ and is minimized for $n \approx \left\lfloor \kappa e \right\rfloor - 1$.
        Choosing such a $n$ leads to a contraction factor upper bounded by $1 - \frac{1}{\kappa e}$.

        This corresponds to the convergence rate obtained by applying $\left\lfloor \kappa e \right\rfloor - 1$ steps of Algorithm~\ref{alg:hb_general} and then restarting it.
        This is described as Algorithm~\ref{alg:hb_restart}.

        $\left\lfloor \kappa e \right\rfloor - 1$ steps of Algorithm~\ref{alg:hb_general} therefore leads to a contraction factor of $\frac{\kappa}{\left\lfloor \kappa e\right\rfloor} \leq \left(1 - \frac{1}{\kappa e}\right)^{\left\lfloor \kappa e \right\rfloor - 1}$.
        Thus applying the same algorithm restarted every $\left\lfloor \kappa e \right\rfloor - 1$ steps leads to a contraction factor of $\left(1 - \frac{1}{\kappa e}\right)^{q \left(\left\lfloor \kappa e \right\rfloor - 1\right)}$ after $q \left(\left\lfloor \kappa e \right\rfloor - 1\right)$ steps.

    $\hfill\blacksquare$

    \begin{center}
        \fbox{\parbox{\textwidth}{
        \begin{restatable}{Cor}{cor:hb_restart}\textbf{\emph{}}
            \label{cor:hb_restart}
            Algorithm~\ref{alg:hb_restart} verifies for every $n$,
            \begin{equation}
                 f(x_n) - f_\star \leq h\left(e \left(1 - \frac{1}{\kappa e}\right)^n d(x_0, \mathcal{X}_\star)^2 \right).
            \end{equation}
        \end{restatable}
        }}
    \end{center}

    \noindent \textit{Proof.}

        Consider $n = q \left(\left\lfloor \kappa e \right\rfloor - 1\right) + r$, with $0 \leq r < \left\lfloor \kappa e \right\rfloor - 1$.

        Combining Theorem~\ref{thm:hb_restart} applied on $q \left(\left\lfloor \kappa e \right\rfloor - 1\right)$ steps and Theorem~\ref{thm:hb_general} applied for the latest $r$ steps from starting point $x_{q \left(\left\lfloor \kappa e \right\rfloor - 1\right)}$, we get

        \begin{eqnarray*}
             f(x_n) - f_\star & \overset{\text{Theorem~\ref{thm:hb_restart}}}{\leq} & h\left(\frac{d(x_{q \left(\left\lfloor \kappa e \right\rfloor - 1\right)}, \mathcal{X}_\star)^2}{r+1} \right) \\
             & \overset{\text{Theorem~\ref{thm:hb_general}}}{\leq} & h\left(\frac{\left(1 - \frac{1}{\kappa e}\right)^{q \left(\left\lfloor \kappa e \right\rfloor - 1\right)} d(x_0, \mathcal{X}_\star)^2}{r+1} \right) \\
             & \leq & h\left(\frac{\left(1 - \frac{1}{\kappa e}\right)^{-r} \left(1 - \frac{1}{\kappa e}\right)^{n} d(x_0, \mathcal{X}_\star)^2}{r+1} \right) \\
             & \leq & h\left(\frac{e \left(1 - \frac{1}{\kappa e}\right)^n d(x_0, \mathcal{X}_\star)^2}{r+1} \right) \\
             & \leq & h\left(e \left(1 - \frac{1}{\kappa e}\right)^n d(x_0, \mathcal{X}_\star)^2 \right)
        \end{eqnarray*}

    $\hfill\blacksquare$

    \begin{Ex}\textbf{\emph{(Logistic regression)}}
        The logistic objective function is strictly convex and smooth.
        But it is not strongly convex.
        However, its square is still convex, and $\QG^+$ (not necessarily smooth anymore), and is still not necessarily strongly convex, but is $\QG^-$.
        Therefore, Theorem~\ref{thm:hb_restart} provides a linear convergence guarantee for Algorithm~$\ref{alg:hb_restart}$.

        Note that without strong convexity, the classical theory does not guarantee linear convergence of Logistic regression without further very specialized analysis.
    \end{Ex}

\end{document}